\documentclass[a4paper,11pt,reqno]{amsart}
\usepackage{graphicx}
\usepackage[english]{babel}

\usepackage{amssymb,enumerate,amsmath}
\usepackage[bookmarks,colorlinks]{hyperref}
\usepackage{tikz}
\usepackage{array}
\usepackage{mathrsfs}
\usepackage[latin1]{inputenc}
\usepackage[T1]{fontenc}
\usepackage{ae,aecompl}

\usepackage[margin=2.5cm]{geometry}

\newcommand{\ZZ}{\mathbb{Z}}
\newcommand{\NN}{\mathbb{N}}
\newcommand{\QQ}{\mathbb{Q}}
\newcommand{\PP}{\mathbb{P}}
\renewcommand{\AA}{\mathbb{A}}
\newcommand{\GrassFunctor}[2]{\underline{\mathbf{Gr}}_{#1}^{#2}}
\newcommand{\GrassScheme}[2]{\mathbf{Gr}_{#1}^{#2}}

\DeclareMathOperator{\reg}{reg}
\DeclareMathOperator{\rk}{rk}
\DeclareMathOperator{\supp}{Supp}

\DeclareMathOperator{\rank}{rk}
\DeclareMathOperator{\Proj}{Proj}
\DeclareMathOperator{\Spec}{Spec}

\newcommand{\xcoeff}{\textnormal{coeff}_{\xx}}

\newcommand{\xx}{x}
\newcommand{\DD}{\Delta}

\newcommand{\Ht}{\textnormal{Ht}}

\newcommand{\PGL}{\textsc{pgl}}

\newcommand{\HilbScheme}[2]{\mathbf{Hilb}^{#2}_{#1}}
\newcommand{\HilbFunctor}[2]{\underline{\mathbf{Hilb}}^{#2}_{#1}}

\newcommand{\Kalg}{(k\textnormal{-Algebras})}
\newcommand{\Sets}{(\textnormal{Sets})}
\newcommand{\schm}{(\textnormal{Schemes}/k)}
\newcommand{\Rings}{(\textnormal{Rings})}

\newcommand{\comp}[1]{{#1}^{c}}
\newcommand{\ee}{\mathrm{e}}
\renewcommand{\a}{\mathrm{a}}
\renewcommand{\b}{\mathrm{b}}

\newcommand{\RevLex}{\mathtt{DegRevLex}}
\newcommand{\cN}{\mathcal N} 
\renewcommand{\geq}{\geqslant}
\renewcommand{\leq}{\leqslant}

\numberwithin{equation}{section}
\newtheorem{theorem}{Theorem}[section]
\newtheorem{corollary}[theorem]{Corollary}
\newtheorem{proposition}[theorem]{Proposition}
\newtheorem{lemma}[theorem]{Lemma}

\theoremstyle{definition}
\newtheorem{definition}[theorem]{Definition}

\newtheorem{remark}[theorem]{Remark}
\newtheorem{notation}[theorem]{Notation}

\begin{document}
\title{Extensors and the Hilbert scheme}

\author[J.~Brachat]{Jerome Brachat}
\address{Jerome Brachat\\INRIA Sophia Antipolis\\ 2004 route des Lucioles, B.P.~93\\ 06902 Sophia Antipolis \\ Cedex France}
\email{\href{mailto:jerome.brachat@inria.fr}{jerome.brachat@inria.fr}}

\author[P.~Lella]{Paolo Lella}
\address{Paolo Lella\\ Universit\`a degli Studi di Torino \\ Dipartimento di Matematica \\ Via Carlo Alberto 10 \\ 10123 Torino \\ Italy}
\email{\href{mailto:paolo.lella@unito.it}{paolo.lella@unito.it}}
\urladdr{\url{http://www.personalweb.unito.it/paolo.lella/}}

\author[B.~Mourrain]{Bernard Mourrain}
\address{Bernard Mourrain\\ INRIA Sophia Antipolis\\ 2004 route des Lucioles, B.P.~93\\ 06902 Sophia Antipolis \\ Cedex France}
\email{\href{mailto:bernard.mourrain@inria.fr}{bernard.mourrain@inria.fr}}
\urladdr{\url{http://www-sop.inria.fr/members/Bernard.Mourrain/}}

\author[M.~Roggero]{Margherita Roggero}
\address{Margherita Roggero\\ Universit\`a degli Studi di Torino \\ Dipartimento di Matematica \\ Via Carlo Alberto 10 \\ 10123 Torino \\ Italy}
\email{\href{mailto:margherita.roggero@unito.it}{margherita.roggero@unito.it}}
\urladdr{\url{http://www.personalweb.unito.it/margherita.roggero/}}


\begin{abstract}  The Hilbert scheme $\HilbScheme{p(t)}{n}$ parametrizes closed subschemes and families of closed subschemes in the projective space
$\PP^n$ with a fixed Hilbert polynomial $p(t)$. It is classically realized as a closed subscheme of a  Grassmannian or a product 
of Grassmannians.   In this paper we consider schemes over a field $k$ of characteristic zero and we present a new proof of the existence of the Hilbert scheme as 
a subscheme of the Grassmannian $\GrassScheme{p(r)}{N(r)}$,  where $N(r)= h^0 (\mathcal{O}_{\PP^n}(r))$. Moreover, we exhibit  explicit equations defining it
in the Pl\"ucker  coordinates of the Pl\"ucker embedding of  $\GrassScheme{p(r)}{N(r)}$.  

Our proof of   existence does not need some of the classical tools used in previous proofs, as flattening stratifications and Gotzmann's Persistence Theorem.
 
The degree of our  equations is $\deg p(t)+2$, lower than  the degree of the equations  given by Iarrobino and Kleiman 
 in 1999 and  also lower (except for the case of hypersurfaces)  than the degree  of those  proved by Haiman 
 and Sturmfels in 2004 after  Bayer's conjecture in 1982. 

The novelty of our approach mainly relies  on  the  deeper attention to the intrinsic symmetries of the Hilbert scheme and on some results about 
Grassmannian based on the notion of extensors.
\end{abstract}

\maketitle

\section*{Introduction}
The study of Hilbert schemes is a very active domain in algebraic geometry. The Hilbert scheme was introduced by Grothendieck \cite{groth} as the scheme
representing the contravariant functor $\HilbFunctor{p(t)}{n}\colon \schm^{\circ}\rightarrow \Sets$ that associates to a scheme $Z$  the set of flat families 
$X \hookrightarrow \PP^n \times_{\Spec k} Z\rightarrow Z$ whose fibers have Hilbert polynomial $p(t)$.
Thus, the Hilbert scheme  $\HilbScheme{p(t)}{n}$  parametrizes the universal family of
subschemes in the projective space $\PP^n$  with  Hilbert polynomial $p(t)$. It is natural to embed the Hilbert functor in a suitable Grassmann
functor and to construct the Hilbert scheme as a
subscheme of a Grassmannian $\GrassScheme{p(r)}{N(r)}$ for a sufficiently large $r$, where $N(t)$ equals $\binom{n+t}{n}$.

Over the years, several authors  addressed the problem  of finding simpler  proofs of the representability of the Hilbert functor  and  
explicit equations for the representing scheme.  This is also the aim of the present paper.  
In fact, we present a new proof of the existence of the Hilbert scheme as a subscheme of 
$\GrassScheme{p(r)}{N(r)}$ and we exhibit explicit equations defining it in the case of a field $k$ of characteristic 0.  

There are some   reasons for which we consider our work significant that concern the tools used in the proofs and the shape of the equations,
in particular the degree. 

In order to simplify Grothendieck's proof, a first crucial point  is the concept of regularity that Mumford introduced for  the  choice of the degree 
$r$    \cite{Mumford,Nitsure}.
 A further simplification is due to Gotzmann, whose Regularity Theorem gives  
 a formula for the  minimum  $r$
only depending on $p(t)$ \cite{Gotzmann}. 
Other key  tools and results  that usually appear in this context are flattening stratifications, fitting ideals,   Gotzmann's Persistence Theorem and 
Macaulay's Estimates on the Growth of Ideals.

In this paper, the number $r$ is always  that given by Gotzmann's formula and in our proof we make use of  Macaulay's Estimates, but we do not need any of 
the other quoted results. We replace them by  a deeper attention to the  inner symmetries of   the  Hilbert scheme   induced by the action  of the   projective
linear group on  $\PP^n$,  and by exploiting the nice combinatorial properties of Borel-fixed ideals.  These are far from being new ideas to study  
Hilbert schemes. Indeed, they play a central role in some of the more celebrated and general achievements on this topic, first of all   Hartshorne's proof of
connectedness \cite{HartThesis}. However, to our knowledge, they have never been used before to prove the existence or to derive equations for $\HilbScheme{p(t)}{n}$.

\medskip

The proof of the representability of the Hilbert functor  given by Haiman and Sturmfels in \cite{HaimSturm}   following Bayer's strategy starts 
with a reduction to   the local case;  the open cover  of 
 $\HilbFunctor{p(t)}{n}$ they  consider  is that  induced by the standard open cover of 
 $\GrassFunctor{p(r)}{N(r)}$. We introduce a new open cover for  the Grassmann   functor,   that we will call the {\it Borel open cover}. It   
 is obtained considering  only a few     open subfunctors $\underline{\mathbf{G}}_{\mathcal{I}} $ of the standard cover,  each  corresponding to 
 a  Borel-fixed ideal  generated by $N(r)-p(r)$ monomials of degree $r$,  and all the open subfunctors $\underline{\mathbf{G}}_{\mathcal{I},g} $, 
 for every $g\in  \PGL(n+1)$,  
 in their  orbit     (Proposition \ref{prop:GrBorelSubFunctors}).  
 The Borel open subfunctors $\underline{\mathbf{H}}_{\mathcal{I},g} $ of the Hilbert functor  are defined accordingly. 

 Restricting to  each Borel subfunctor $\underline{\mathbf{G}}_{\mathcal{I},g} $, the properties of   $J$-marked sets and bases  over  a Borel-fixed ideal  $J$ developed in \cite{LR2} allow us to  prove that  
 $\underline{\mathbf{H}}_{\mathcal{I},d} $
 is representable and to obtain a new proof of  the existence of  $\HilbScheme{p(t)}{n}$   (Theorem \ref{representable}).

\medskip

Towards the aim of deriving equations for the Hilbert scheme, we  then expand  the notion of marked set  to the universal element of the family 
 \[
 \mathcal{F} \hookrightarrow \PP^n \times_{\Spec k} \GrassScheme{p(r)}{N(r)} \rightarrow \GrassScheme{p(r)}{N(r)}
 \]
 parameterized by the 
Grassmannian  and to its exterior powers. Indeed,   exploiting  the  notion of an \emph{extensor} and its properties, we  
obtain a description of the universal element   
by a  set of bi-homogeneous generators  of bi-degree $(r,1)$ in $k[\xx,\DD]$, where $\xx$ and  $\DD$  are compact notation for
the set of variables on $\PP^n$ and the Pl\"ucker coordinates on the Grassmannian. We also obtain a 
similar description  (again linear w.r.t.~$\DD$)  for 
 the    exterior powers of the universal element of $\mathcal{F} $ (Theorem \ref{universal}). These sets of generators  allow us to  write explicitly
 a set of  equations 
 for the Hilbert scheme in the ring $k[\DD]$ of the Pl\"ucker coordinates (Theorem \ref{th:mainTheorem}). 

 \medskip
 
The degree of our equations is upper bounded by $d+2$, where $d:=\deg p(t)$ is the dimension of the subschemes of $\PP^n$ parametrized by $\HilbScheme{p(t)}{n}$. 
It is interesting that the degree of the equations is so close to the geometry of the involved objects. Furthermore, 
 $d+2$ is  lower than the degree of the  other known sets of equations for the embedding 
of the Hilbert scheme in a single Grassmannian. 
We quote  the  
equations of degree $N(r+1)-p(r+1)+1$  in local coordinates given by  Iarrobino and Kleiman \cite[Proposition C.30]{IarrobinoKleiman}, 
and the  equations of degree $n+1$ in the Pl\"ucker
 coordinates  conjectured by Bayer in 
his thesis 
 \cite{Bayer82} and obtained by  Haiman and
 Sturmfels  as a special case of a more general result in \cite{HaimSturm}.   
 
 By the way, we observe that our method, applied with slightly different strategies, also allows to  obtain sets of 
 equations very similar to those by  Iarrobino and Kleiman and by  Haiman and
 Sturmfels (Theorems \ref{th:IarroKleiman} and \ref{th:BayerHaimanSturmfels}).    

At the end of the paper we apply our results in order to compute a set of equations defining the Hilbert schemes of 2 points in $\PP^2$, $\PP^3$ and $\PP^4$ and of 3 points in $\PP^2$. In particular, we illustrate in detail our method in the case of $\HilbScheme{2}{2}$ and we compare the equations we obtain with those obtained by Brodsky and Sturmfels \cite{BrodskySturmfels}. We observe that the two sets of equations, though different, generate
 the same ideal, more precisely the saturated ideal of $\HilbScheme{2}{2}$ in $\GrassScheme{2}{6} \subset \PP^{14}$. Our equations describe the saturated ideal also in the case of $\HilbScheme{2}{3}$ in $\GrassScheme{2}{10} \subset \PP^{44}$, $\HilbScheme{2}{4}$ in $\GrassScheme{2}{15} \subset \PP^{104}$ and $\HilbScheme{3}{2}$ in $\GrassScheme{3}{10} \subset \PP^{119}$, but we do not know if this nice property holds in general. However, the lower degree marks a significant step
forward in order to compute this special ideal (see Table \ref{tab:compare}) and allows further experiments and investigations.  

 \medskip

 Let us now explain the structure of the paper. 
  In Section \ref{sec:HilFunctor}, we 
 recall some properties 
 that we will use throughout the
 paper. In particular,  we describe the Hilbert functor, its relation with the Grassmann functor and  the standard open cover.
In Section \ref{sec:Borel},  we introduce the Borel open cover.  
 Section \ref{sec:representability} contains  the generalities about  marked sets and bases over Borel-fixed ideals and it ends with the first main result
 of the paper, namely   Theorem \ref{representable} on the representability of the Hilbert functor.
Section \ref {sec:extensors} contains the results on Grassmannians based  on the theory of extensors (Theorem \ref{universal}).  
In Section \ref{sec:BLMR}, after  some new technical results about marked bases,     we present the equations  defining the Hilbert scheme and we prove  
their correctness (Theorem  \ref{th:mainTheorem}). In Subsections  \ref{subsec:IK}  and \ref{subsec:HS}  we  derive 
equations similar to those   by Iarrobino-Kleiman and by Haiman-Sturmfels. 
In Section \ref{sec:example}, we illustrate  the constructions and results of the paper in the case of Hilbert schemes describing 2 or 3 points.


\section{Notation}\label{sec:Gen}

 Let  $k$   be   a field of characteristic 0. In the following $k[\xx]$ will denote the polynomial ring $ k[x_0,\ldots,x_n]$ and  $\PP^n$ 
 the $n$-dimensional projective space $\Proj k[\xx]$. For a  
$k$-algebra $A$, we will denote by $A[\xx] := A\otimes_k k[\xx]$ the polynomial ring with coefficients in $A$ and by $\PP^n_A$ 
the projective space $\Proj A[\xx]=\PP^n\times_{\Spec k}\Spec A$.  As usual, for a subset $E$ of a ring $R$, we denote by  $(E)$ the ideal 
of $R$ generated by $E$ and  for a subset $F$ of an $R$-module $M$, we denote by  $\langle F \rangle$ the $R$-submodule of $M$ generated by $E$; 
we sometimes  write ${}_R(E )$ and  ${}_R\langle F \rangle$   when more than one ring is involved.

Let us now consider a scheme $X \subset \PP^n_A$. For each prime ideal $\mathfrak{p}$ of $A$, we denote by $A_{\mathfrak{p}}$ the 
localization
 in $\mathfrak{p}$, by $k(\mathfrak{p})$ the residue field and by $X_{\mathfrak{p}}$ the fiber of the structure morphism
 $X \rightarrow \Spec A$.
 The Hilbert polynomial $p_{\mathfrak{p}}(t)$ of $X_{\mathfrak{p}}$ is defined as
\[
p_{\mathfrak{p}}(t) = \dim_{k(\mathfrak{p})} H^0 (X_{\mathfrak{p}},\mathcal{O}_{X_{\mathfrak{p}}}(t)) \otimes_{k} k(\mathfrak{p}),
\qquad t \gg 0.
\]
 If $X$ is flat over $\Spec A$ and  the Hilbert polynomial $p_{\mathfrak{p}}(t)$ of every localization coincides with $p(t)$,
then $p(t)$ is 
called the Hilbert polynomial of $X$ 
(for further details see \cite[III, \S 9]{Hartshorne}).  There exists a positive integer $r$ only depending on $p(t)$,
called \emph{Gotzmann number}, for which the ideal sheaf $\mathcal{I}_X$ of each scheme $X $ with Hilbert polynomial $p(t)$ is $r$-regular
(in the sense of 
Castelnuovo-Mumford regularity). 
By Gotzmann's Regularity Theorem (\cite[Satz (2.9)]{Gotzmann} and \cite[Lemma C.23]{IarrobinoKleiman}), this implies   the surjectivity of the morphism
\[
H^0\big(\mathcal{O}_{\PP^n_A}(r)\big)\ \xrightarrow{\phi_X}\ H^0\big(\mathcal{O}_X(r)\big).
\]

We will denote by $N(t)$ the dimension of  $k[\xx]_t$. 
The polynomial $q(t) := N(t) - p(t)$ is  the Hilbert polynomial of the saturated ideal
defining $X$ and it is called the \emph{volume polynomial} of $X$. In particular, for $t=r$ the Gotzmann number of $p(t)$, we set
$p:=p(r)$, $q:=q(r)$ and $N:=N(r)$.

We will use the usual notation for terms  $x^\alpha := x_0^{\alpha_0} \cdots
 x_n^{\alpha_n}$, where $\alpha = (\alpha_0,\ldots,\alpha_n) \in \NN^{n+1}$. When a term  order  comes into play, we  assume the variables 
 ordered as  $x_0 <  \cdots  < x_n$;   we will denote by $<_{\RevLex}$ and 
$<_{\mathtt{Lex}}$ the 
degree reverse lexicographic and the lexicographic orders. We will denote
by $x^{\alpha(i)}$ the $i$-terms of degree $r$ in descending $\RevLex$ order. 
For any term $x^\alpha$,  let $\min(x^\alpha)$ and $\max(x^\alpha)$ denote respectively the minimal and  the maximal variable
which divides $x^\alpha$.

 For any polynomial $f\in A[\xx]$,  
the support 
$\supp(f)$ of $f$ is the set of terms  that appear in $f$ with  non-zero coefficient and $\xcoeff(f)\subset A$ is the set  of coefficients 
of the  terms in $\supp(f)$; with the obvious meaning, we use  the notation $\xcoeff(U)$ also if $U$ is  a subset of   $ A[\xx]$. 

We loosely denote by the same letter the monomial ideals in $k[\xx]$ and that in $A[\xx]$ generated by the same set of terms.  If $J$ 
is a monomial ideal, we will denote by $B_J$ its minimal monomial basis and by $\cN(J)$ the set of terms in $k[\xx]\setminus J$. For a subset $V$ of a standard
graded module $R=\bigoplus_{t}R_t$,   $V_s$ and $V_{\geq s}$  will denote respectively $V\cap R_s$ and $V\cap \bigoplus_{t\geq s}R_t$.

An $s$-multi-index  $\mathcal{H}=(h_1,\ldots,h_s) $  is an ordered sequence   $h_1 < h_2 <\dots <h_s$ in $\{1, \dots, N\}$; 
its complementary  $\comp{\mathcal{H}}$ is the 
$(N-s)$-multi-index  with entries in the set  $\{1,\ldots,N\} \setminus \mathcal{H}$.  For any $s$-multi-index $\mathcal{H}$, 
 we will denote by $\varepsilon_{\mathcal{H}} \in 
 \{-1,1\}$ the signature of the permutation $(1,\ldots,N) \mapsto {\mathcal{H}}{,}\comp{\mathcal{H}}$.
Moreover,   if $\mathcal H\subset \mathcal K$,  
$\varepsilon_{\mathcal H}^{\mathcal K}$ is  the signature of $\mathcal K \mapsto \mathcal H,\mathcal K \setminus \mathcal H$.
  For every $m\leq N-p$, we  will denote by $\mathcal{E}^{(m)} $  the set of all $(p+m)$-multi-indices.

For every $\mathcal I\in \mathcal E^{(0)}$, $J(\mathcal I)$ is the ideal generated by the terms $x^{\alpha(j)}$ corresponding to the indices
$j\in  \comp{\mathcal I}$.


\section{Hilbert and Grassmann functors}\label{sec:HilFunctor}

In the following, $\HilbFunctor{p(t)}{n}$ will denote the Hilbert functor $\schm^{\circ}\rightarrow \Sets$
that associates to an object $Z$ of the category of schemes over $k$ the set
\begin{equation*}
 \HilbFunctor{p(t)}{n}(Z) = \{X \subset \PP^n \times_{\Spec k} Z\ \vert\ X \rightarrow Z \text{ flat  with  Hilbert
 polynomial } p(t)\}.
\end{equation*}
and to any morphism of schemes $f\colon Z \rightarrow Z'$ the map 
\[
\begin{split}
\HilbFunctor{p(t)}{n}(f)\colon&\ \HilbFunctor{p(t)}{n}(Z') \rightarrow \HilbFunctor{p(t)}{n}(Z)\\
&\parbox{2.18cm}{\centering $X'$} \mapsto\ X' \times_{Z'} Z
\end{split}
\]

It is easy to prove that $\HilbFunctor{p(t)}{n}$ is a Zariski sheaf \cite[Section 5.1.3]{Nitsure}; 
hence, we can consider it as a covariant functor from the category
of noetherian  $k$-algebras \cite[Lemma E.11]{Sernesi}
\begin{equation*}
\HilbFunctor{p(t)}{n}\colon \Kalg \rightarrow \Sets
\end{equation*}
such that for every finitely generated  $k$-algebra $A$
\[
\HilbFunctor{p(t)}{n}(A) = \left\{ X \subset \PP^n_A\ \vert\ X \rightarrow \Spec A \text{ flat  with  Hilbert 
polynomial }
 p(t) \right\}.
\]
and for any $k$-algebra  morphism $f\colon A \rightarrow B$
\[
\begin{split}
\HilbFunctor{p(t)}{n}(f)\colon&\ \HilbFunctor{p(t)}{n}(A)\ \rightarrow \parbox{3cm}{\centering $\HilbFunctor{p(t)}{n}(B)$}\\
& \parbox{2.2cm}{\centering $X$} \mapsto\ X \times_{\Spec A} \Spec B.
\end{split}
\]

The Hilbert scheme $\HilbScheme{p(t)}{n}$ is defined as the scheme representing the Hilbert
functor. Our notation for the Hilbert functor follows   that  used for instance in 
\cite{HaimSturm}, 
 where the functor of points of a scheme $Z$ is denoted by $\underline{Z}$. Note that    we are not assuming the representability of    
 $\HilbFunctor{p(t)}{n}$ as a known fact, but we will prove it  at the end of Section \ref{sec:representability}.

 Let us briefly recall the strategy of the construction of the Hilbert scheme based on Castelnuo-vo-Mumford regularity and Gotzmann number.
The following proposition suggests to look for an embedding in a representable functor and reduce to the local case.

\begin{proposition}[{\cite[Proposition 2.7 and Corollary 2.8]{HaimSturm}}]\label{prop:openSubfunctors}
Let $Z$ be a scheme and ${\eta}\colon {\mathbf{F}} \rightarrow \underline{Z}$ be a natural transformation of functors $\Kalg \rightarrow 
\Sets$, where 
 ${\mathbf{F}}$ is a Zariski sheaf. Suppose that  $\underline{Z}$ has a cover of open subsets $\underline{U}$ such that each subfunctor 
 ${\eta}^{-1}(\underline{U}) 
 \subseteq 
{\mathbf{F}}$ is representable. Then, also  ${\mathbf{F}}$ is representable.

 Moreover, if the natural transformations ${\eta}^{-1}(\underline{U}) \rightarrow 
\underline{U}$, given by restricting ${\eta}$, are induced by closed embeddings of schemes, then so is ${\eta}$.
\end{proposition}

The overall strategy introduced by Bayer \cite{Bayer82} for the construction of the Hilbert scheme uses an embedding in a  Grassmann 
functor (for a detailed discussion we refer Section 2 of \cite{HaimSturm}  and to \cite[Section VI.1]{EisenbudHarris}).
If  $X \in \HilbFunctor{p(t)}{n}(A) $, then  by flatness $H^0(\mathcal{O}_X(r))$ is a locally free $A$-module of rank $p(r)$. Hence,  
the surjective  map 
 $\phi_X \colon H^0\big(\mathcal{O}_{\PP^n_A}(r)\big) \simeq A^N \rightarrow\ H^0\big(\mathcal{O}_X(r)\big) $   is an element of the set defined by 
the Grassmann 
functor $\GrassFunctor{p}{N}$ over $A$. Indeed, the Grassmann functor 
$
\GrassFunctor{p}{N}\colon \Kalg  \rightarrow \Sets
$
associates to every finitely generated $k$-algebra  $A$ the set
\begin{equation*}
\begin{split}
 \GrassFunctor{p}{N}(A)=\left\{\begin{array}{c}\text{isomorphism classes of epimorphisms} \\  \pi \colon A^N \rightarrow P\text{ of
 locally free modules of rank }
p\end{array} \right\} .
\end{split}
\end{equation*}
Equivalently, we can define
\begin{equation}\label{eq:sottomodulo}
\begin{split}
 \GrassFunctor{p}{N}(A) =\left\{\begin{array}{c}\text{submodules  } L\subseteq A^N  \text { such that} \\ A^N/L  \text{ is locally free of rank }
p\end{array} \right\}.
\end{split}
\end{equation}
In the second formulation, $\pi$ is the canonical projection $\pi_L\colon A^N\rightarrow A^N/L$.
This functor is representable and the representing scheme $\GrassScheme{p}{N}$ is called the Grassmannian (see \cite[Section 16.7]{Vakil}).

We  fix the canonical basis  $\{\a_1,\ldots,\a_N\}$  for $A^N$ and the isomorphism   $A^N\simeq  H^0\big(\mathcal{O}_{\PP^n_A}(r)\big)$ 
given by  $\a_i \mapsto x^{\alpha(i)} $. Thus,  we obtain  a universal family 
\begin{equation}\label{eq:universalfamily}
\mathcal{F} \hookrightarrow \PP^n \times \GrassScheme{p}{N} \rightarrow \GrassScheme{p}{N}
\end{equation} 
parameterized by the 
Grassmannian and the natural transformation of functors 
\begin{equation*}
\mathscr{H}\colon  \HilbFunctor{n}{p(t)}\ \rightarrow\ \GrassFunctor{p}{N}
\end{equation*}
sending
 $X \in \HilbFunctor{n}{p(t)}(A)$ to $\pi_X \colon H^0\big(\mathcal{O}_{\PP^n_A}(r)\big) \rightarrow H^0\big(\mathcal{O}_X(r)\big) \in \GrassFunctor{p}{N}(A)$ (or equivalently   to $ L=H^0\big(\mathscr{I}_X(r)\big)$).

The Grassmannian has the following well-known open cover that we call the \emph{standard open cover} of $\GrassFunctor{p}{N}$. 
Let us fix a basis  $\{\ee_1,\ldots,\ee_p\}$  for $A^p$.   For every $\mathcal{I}=(i_1, \dots, i_p)
\in \mathcal E^{(0)}$, 
let us consider the injective morphism 
\[
\begin{split}
\Gamma_{\mathcal{I}}\colon A^p &\rightarrow A^N \\
\ee_j & \mapsto \a_{i_j}
\end{split}
\]
and the  subfunctor $\underline{\mathbf{G}}_{\mathcal{I}}$ that associates to every noetherian $k$-algebra $A$ the set 
\begin{equation*}
\underline{\mathbf{G}}_{\mathcal{I}}(A) = \left\{ L\in \GrassFunctor{p}{N}(A) \text{ such that } \pi_L \circ \Gamma_{\mathcal{I}} 
\text{ is surjective}\right\}.
\end{equation*}

 \begin{proposition}\label{prop:openCoverGrass} For  ${\mathcal I}\in {\mathcal E}^{(0)}$, the $ \underline{\mathbf{G}}_{\mathcal{I}}$ 
 are open subfunctors of  $\GrassFunctor{p}{N}$ that  cover it.
\end{proposition}
\begin{proof}     See \cite[Section 22.22]{stacks-project}.
\end{proof}

\begin{remark}\label{rk:precisazioni} For every  $L \in \underline{\mathbf{G}}_{\mathcal{I}}(A)$ the map $\pi_L \circ \Gamma_{\mathcal{I}} $
is an isomorphism, as it is a surjective morphism from a free $A$-module to a locally free $A$-module of the same rank. Therefore, $L$ is 
the kernel of the epimorphism $\phi_L:=(\pi_L \circ \Gamma_{\mathcal{I}})^{-1}\circ \pi_L \colon A^N \rightarrow A^p$ such that $\phi_L(\a_{i_j} )=
\ee_j$ for every ${i_j} \in \mathcal I$. 

On the other hand, the kernel of every  surjective morphisms $\phi\colon A^N \rightarrow A^p$ sending $\a_{i_j}$ to $\ee_j$ is by definition a
module $L\in \underline{\mathbf{G}}_{\mathcal{I}}(A)$: we will write $\phi_L$ instead of $\phi$ to emphasize this correspondence.

Every map $\phi_L$ is  completely determined  
 by the images of the   $q=N-p$ elements   $\a_h$ with $h \in \comp{\mathcal{I}}$. 
 If $\phi_L(\a_h) =\sum_{j=1}^p \gamma_{hj} \, \ee_j= \phi_L\left(\sum_{j=1}^p \gamma_{hj}\, \a_{i_j}\right)$,   the kernel $L$ 
 contains the free $A$-module $L'$
 generated by the $q$ elements 
$\b_h:=\a_h -\sum_{j=1}^p \gamma_{hj}\, \a_{i_j}\in A^N$ with $h \in \comp{\mathcal{I}}$. Then,   $A^N= L'\oplus \langle \a_{i_j} \ 
\vert \ i_j\in \mathcal I \rangle\subseteq  L\oplus \langle \a_{i_j} \ \vert \ i_j\in \mathcal I \rangle\subseteq A^N$, so that    $L=L'$ 
and  $A^N/L$ are  free $A$-modules of rank $q$ and $p$ respectively.

Through  the fixed isomorphism $A^N\simeq  H^0\big(\mathcal{O}_{\PP^n_A}(r))$ given by $\a_j\mapsto x^{\alpha(j)}$, the elements $\b_h$ 
correspond to  polynomials $ f_{\alpha(h)}:=x^{\alpha(h)} -\sum_{j=1}^p \gamma_{hj}x^{\alpha(i_j)} \in H^0\big(\mathcal{O}_{\PP^n_A}(r)\big)$. 
In this way, for  $L={\mathscr{H}}(A)(X)\in  \underline{\mathbf{G}}_{\mathcal{I}}(A)$,  the polynomials $f_{\alpha(h)}$ generate  the 
ideal  $(I_X)_{\geq r}$, while   for a general $L\in  \underline{\mathbf{G}}_{\mathcal{I}}(A)$, the $A$-module
$ \langle f_{\alpha(h)}, h \in \comp{\mathcal{I}}\rangle \subseteq H^0\big(\mathcal{O}_{\PP^n_A}(r)\big)$  is    free   of rank $q$,
but the Hilbert polynomial of  $\Proj(A[\xx]/( f_{\alpha(h)}, h \in \comp{\mathcal{I}}))$ is not necessarily $p(t)$. 

 In the following, keeping in mind the above construction, we often consider the ideal $I =( f_{\alpha(h)}, h \in \comp{\mathcal{I}})$ as an
 element of  $\underline{\mathbf{G}}_{\mathcal{I}}(A)$, identifying it with the $A$-module $L=I_r$. In the same way, we will
 write $I\in \underline{\mathbf{H}}_{\mathcal{I}}(A)$ when   $I\in \underline{\mathbf{G}}_{\mathcal{I}}(A)$ and  the Hilbert polynomial  $\Proj( A[\xx]/I)$
 is $p(t)$.
 \end{remark}

The proof of the representability of the Hilbert functor  after Bayer's strategy given in \cite{HaimSturm} uses the  open cover of 
$\HilbFunctor{p(t)}{n}$ of  the subfunctors  
$\underline{\mathbf{H}}_{\mathcal{I}}:={\mathscr{H}}^{-1}(\underline{\mathbf{G}}_{\mathcal{I}})
\cap \HilbFunctor{p(t)}{n}$, that we will call the {\it standard open cover} of $\HilbFunctor{p(t)}{n}$.
  
In this paper  we introduce new  open covers of the Grassmann and the Hilbert functors, called \emph{Borel open covers}, that  
take into account of the action of the projective linear group on the Grassmann and Hilbert functors induced by that on $\PP^n$.


\section{The Borel open cover}\label{sec:Borel}

An ideal $J \subset k[\xx]$ is said \emph{Borel-fixed} if it is fixed by the action of the Borel subgroup of the upper triangular matrices. 

These ideals are  involved in many  general results about Hilbert schemes for the following reason. Galligo \cite{galligo} and Bayer and Stillman
\cite{BayerStillmanGIN} proved that the generic initial ideal of any ideal is Borel-fixed, which means, in the context of Hilbert schemes, 
that any component and any intersection of components of $\HilbScheme{p(t)}{n}$ contains at least a point corresponding to a scheme defined by a Borel-fixed ideal. 

In characteristic zero, the notion of Borel-fixed ideals coincide with the notion of strongly stable ideals. An ideal $J$ is said \emph{strongly stable} if, 
and only if, it is generated by terms and
for each term $x^\alpha \in J$  also the term $\frac{x_j}{x_i} x^\alpha$ is 
in  $J$ for all $x_i \mid x^\alpha$ and $x_j > x_i$.   Moreover, the regularity of $J$ is equal to the maximum degree of 
terms in its minimal monomial basis \cite[Proposition 2.11]{Green}. For further details about Borel-fixed ideals see \cite{BLR,Green,MS}.

\begin{notation} 
For any Hilbert polynomial $p(t)$ and for  the related integers $r$, $p$, $N$, $q$ 
\begin{itemize}
\item $\mathbb{B}$ is the set of  the  Borel-fixed ideals in  $k[\xx]$  generated by $q$ terms of degree $r$.
\item   $\mathbb{B}_{p(t)}$ is  the set of  Borel-fixed ideals in $\mathbb{B}$  with Hilbert polynomial $p(t)$.
\item A \emph{Borel multi-index} $\mathcal{I}$ is any multi-index in $ \mathcal E^{(0)}$ such that  $J(\mathcal I)\in \mathbb{B}$.
\item  For every element $g\in \PGL:=\PGL_\QQ(n+1)$, $\widetilde{g}$ denotes the automorphism induced by $g$ on $A[\xx]_r$ and on the Grassmann and
Hilbert functors and $g\centerdot$ denotes  the
corresponding action on an element. 
\end{itemize}
\end{notation}

Notice that the set of Borel-fixed ideals in $\mathbb{B}_{p(t)}$ can be efficiently computed by means of the algorithm presented in \cite{CLMR}
and subsequently improved in \cite{LellaBorel}.
\medskip

For any $p$-multi-index 
$\mathcal{I}\in \mathcal E^{(0)}$ and any $g \in \PGL$, we  consider the following subfunctor  of the Grassmann functor :
\[
\underline{\mathbf{G}}_{\mathcal{I},g}(A) = \left\{ \begin{array}{c} \text{free quotient } A^{N} 
\xrightarrow{\pi_L} A^{N}/L\text{ of rank } p \\ \text{such that } \pi_L \circ \widetilde{g}\circ \Gamma_{\mathcal{I}}
 \text{ is surjective}\end{array} \right\}.
\]

These subfunctors are open, because the functorial automorphism of $\GrassFunctor{p}{N}$ induced by $\widetilde{g}$ extends to
 $\underline{\mathbf{G}}_{\mathcal{I},g} \simeq \underline{\mathbf{G}}_{\mathcal{I},\mathrm{id}} = \underline{\mathbf{G}}_{\mathcal{I}}$. 
 It is obvious that these subfunctors also cover $\GrassFunctor{p}{N}$, but in fact it is sufficient to consider a smaller subset.

\begin{proposition}\label{prop:GrBorelSubFunctors}
The collection of subfunctors
\begin{equation*}
\left\{\underline{\mathbf{G}}_{\mathcal{I},g} \ \big\vert\ g \in \PGL,\ 
 \mathcal{I}\in \mathcal E^{(0)} \text{ s.t.~}J(\mathcal I) \in \mathbb{B} \right\}
\end{equation*}
covers the Grassmann functor $\GrassFunctor{p}{N}$ and 
the representing schemes ${\mathbf{G}}_{\mathcal{I},g}$ cover the Grassmannian $\GrassScheme{p}{N}$.
\end{proposition}
\begin{proof}
Let $\pi \colon A^{N} \rightarrow P$ be an element of $\GrassFunctor{p}{N}(A)$. Following \cite[Lemma 22.22.1]{stacks-project}, 
we prove the result showing that for any $\mathfrak{p}
 \in \Spec A$ there exist a multi-index $\mathcal{I}$ and a change of coordinates $g$ such that the morphism $\pi \circ \widetilde{g}\circ 
\Gamma_{\mathcal{I}}$ is surjective in a neighborhood of $\mathfrak{p}$.

Let $A_{\mathfrak{p}}$ be the local algebra obtained by localizing in $\mathfrak{p}$, $\mathfrak{m}_{\mathfrak{p}}$ its maximal ideal and 
$k(\mathfrak{p})$ the residue field. Tensoring by $k(\mathfrak{p})$ the morphism $\pi$, we obtain the morphism of vector spaces
\[
\pi_{\mathfrak{p}}\colon k(\mathfrak{p})^{N} \rightarrow P_{\mathfrak{p}}/\mathfrak{m}_{\mathfrak{p}}P_{\mathfrak{p}}
\]
whose kernel is a vector subspace of $k(\mathfrak{p})\otimes S_r$ of dimension $q$. 

Now, consider the ideal $I \subset k(\mathfrak{p})\otimes S$ generated by $\ker \pi_{\mathfrak{p}}$ and let $J$ be its generic initial ideal.
 We fix an element $g \in \PGL$ such that $J = \textnormal{in}(g\centerdot I)$. By properties of Gr\"obner
 bases, we know that $\dim_{k(\mathfrak{p})} J_r = \dim_{k(\mathfrak{p})} (g\centerdot I)_r$ ($J$ and $g\centerdot I$ have the same Hilbert 
 function).
Furthermore, the terms of degree $r$ not belonging to $J$ are a basis both of $(k(\mathfrak{p})\otimes S_r)/J_r$ and 
$(k(\mathfrak{p})\otimes S_r)/(g\centerdot I)_r$.

Finally, the multi-index $\mathcal{I}$ is the one such that $J(\mathcal I)=J$. 
\end{proof}

 \begin{definition}\label{BorelcoverGrass}
We call \emph{Borel subfunctor} of
 $\GrassFunctor{p}{N}$
 any element of the collection of subfunctors of Proposition \ref{prop:GrBorelSubFunctors}. 
Moreover, we denote by $\underline{\mathbf{H}}_{\mathcal{I},g}$ the open subfunctor ${\mathscr{H}}^{-1}
 (\underline{\mathbf{G}}_{\mathcal{I},g})
 \cap \HilbFunctor{p(t)}{n}$.
\end{definition}

\begin{theorem}\label{lem:HilbBorelSubFunctors}  
The collection of subfunctors
\begin{equation}\label{eq:HilbBorelSubFunctors}
\left\{\underline{\mathbf{H}}_{\mathcal{I},g} \ \big\vert\ g \in \PGL,\  \mathcal{I}\in \mathcal E^{(0)} \textit{ s.t.~}J(\mathcal I) 
\in \mathbb{B}_{p(t)}
\right\}
\end{equation}
covers the Hilbert functor $\HilbFunctor{p(t)}{n}$. 
\end{theorem}
\begin{proof}
Consider an element $X \in \HilbFunctor{p(t)}{n}(A)$. As above, it is sufficient to prove that for any $\mathfrak{p} \in \Spec A$, there exists a 
subfunctor $\underline{\mathbf{H}}_{\mathcal{I},g}$ such that $X_{\mathfrak{p}} = X \times_k \Spec k(\mathfrak{p})$ is an element of 
$\underline{\mathbf{H}}_{\mathcal{I},g}\big(k(\mathfrak{p})\big)$.

Localizing at $\mathfrak{p}$, we obtain a scheme $X_{\mathfrak{p}}$ flat over $\Spec k(\mathfrak{p})$ with Hilbert polynomial $p(t)$,
 as the flatness and so the Hilbert polynomial are preserved by localization. Let $I_X \subset k(\mathfrak{p})\otimes S$ be the saturated 
ideal defining
 $X_{\mathfrak{p}}$, $I:=(I_X )_{\geq r}$ and $J$ the generic initial ideal of $I$. By the same argument used in the proof of
 Proposition \ref{prop:GrBorelSubFunctors}, 
we fix a change of coordinated $g \in \PGL$ such that $J = \textnormal{in}(g\centerdot I)$ and the 
multi-index $\mathcal{I}\in \mathcal E^{(0)}$ such that $J(\mathcal I)=J$. By construction, $J(\mathcal I) \in \mathbb{B}_{p(t)}$, as $J$ and $I$ 
share the same Hilbert function. 
\end{proof}

\begin{definition}\label{BorelCoverHilb}
The 
 \emph{Borel cover} of
 $\HilbFunctor{p(t)}{n}$ is  the collection of the open subfunctors \eqref{eq:HilbBorelSubFunctors} of Theorem \ref{lem:HilbBorelSubFunctors}. 
\end{definition}

In next section we will prove that  the open subfunctor  $\underline{\mathbf{H}}_{\mathcal{I},g}$ is empty if  $J(\mathcal I) \in \mathbb{B}
\setminus \mathbb{B}_{p(t)}$.


\section{Representability}\label{sec:representability}

Our proof that the Hilbert functor  is representable mainly uses the theory of marked sets and  bases on a Borel-fixed ideal developed in \cite{CR,BCLR,LR2}. 
We recall some of the results and notation contained
in the quoted papers.

\begin{definition}[{\cite[Definitions 1.3, 1.4]{CR}}]
A \emph{monic marked polynomial} (marked polynomial for short) is a polynomial $f\in A[\xx]$ together with a specified term
$x^\alpha$ of $\supp(f)$,  called \emph{head term}
of $f$ and denoted by $\Ht(f)$.  We assume furthermore that the coefficient of $x^\alpha$ in $f$ is $1_A$. Hence, we can write a marked polynomial 
as $f_\alpha = 
x^\alpha-\sum c_{\alpha\gamma}x^\gamma$, with $x^\alpha =\Ht (f_\alpha)$, $x^\gamma \neq x^\alpha$ and $c_{\alpha\gamma}\in A$. 
\end{definition}

\begin{definition}\label{def:markedset}
Let $J$ be a monomial ideal. A finite set $F$ of homogeneous marked polynomials $f_\alpha = x^\alpha-\sum c_{\alpha\gamma}x^\gamma$, 
with $\Ht(f_\alpha) =x^\alpha$, 
is called a \emph{$J$-marked set} if the head terms $x^\alpha$ form the  minimal monomial basis $B_J$ of $J$, 
 and every $x^\gamma$ is an element of  $\cN(J)$. Hence,  $\mathcal N(J)$  generates the quotient $A[\xx]/(F)$ as an $A$-module.

 A $J$-marked set $F$ is a \emph{$J$-marked basis} if  the quotient $A[\xx]/(F)$ is freely generated by $\mathcal N(J)$ as an $A$-module,
 i.e.~$A[\xx]={}_{A[\xx]}(F)\oplus {}_A\langle \mathcal N(J)\rangle$.
\end{definition}

\begin{remark} Observe that if $I$ is generated by a    $J$-marked basis, then   $\Proj(A[\xx]/I)$ is $A$-flat,   since  $A[\xx]/I$ is a 
free $A$-module.

In the following we will consider only  $J$-marked sets $F$  with $J\in \mathbb B$, i.e.~of the shape
\begin{equation} \label{eq:markedset}
F=\left\{  
f_\alpha:= x^\alpha-\sum c_{\alpha\gamma}x^\gamma \ \big\vert\ x^\alpha \in J_r, \ x^\gamma \in \cN(J)_r , \  c_{\alpha\gamma} \in A \right\}.
\end{equation}
For every ideal $I$   generated by such a $J$-marked set $F$, we have  
$A[\xx]_r=\langle {F}\rangle \oplus \langle \cN(J)_r\rangle$, hence  $I_r$ is a free direct summand  of rank    $q$ 
of $A[\xx]_r=H^0 \big(\mathcal{O}_{\PP^n(A)}(r)\big)$ 
and  it corresponds to an element of $\GrassFunctor{p}{N}(A) $. In fact, if $\mathcal I\in \mathcal {E}^{(0)}$ is the 
$p$-multi-index such that $J(\mathcal I)=J$, 
then $I \in \underline{\mathbf{G}}_{\mathcal{I}}(A)$.

 Moreover   $I \in \underline{\mathbf{H}}_{\mathcal{I}}(A)$ if, and only if, 
the Hilbert polynomial of $A[\xx]/(I)$ 
is $p(t)$. Now we will prove that this happens if, and only if, $J\in \mathbb B_{p(t)}$ and $F$ is a $J$-marked basis.
\end{remark}

We need some more properties  concerning Borel-fixed ideals and marked bases.

\begin{lemma}\label{lemma:monomialsBorel} If  $J \in \mathbb{B}$, then  $\rk (J_t) \geq q(t)$ and $k[x_{d+1},\ldots,x_n]_t \subset J_t$
for all $ t \geq r$.
\end{lemma}
\begin{proof} 
The first assertion follows by
Macaulay's Estimates on the Growth of Ideals \cite[Theorem 3.3]{Green}. Thus, the degree of the Hilbert polynomial of $A[\xx]/J$ is at most $d=\deg p(t)$.
By \cite[Proposition 2.3]{CLMR}, we have $x_{d+1}^{r} \in J$ that implies $k[x_{d+1}, \dots, x_n]_t \subset J_{t}$ for $t\geq r$, 
by the strongly stable property.
\end{proof}

\begin{proposition}[{\cite[Lemma 1.1]{EliahouKervaire}, \cite[Lemma 1.2]{BCLR}}]
Let $J$ be a strongly stable ideal and let $B_J$ be its minimal monomial basis.
\begin{enumerate}[(i)]
\item\label{it:borelProducts_i} Each term $x^\alpha$ can be written uniquely as a product $x^{\gamma}x^\delta$ with $x^\gamma \in B_J$ and
$\min x^\gamma
\geq \max x^\delta$. Hence, $x^\delta <_\mathtt{Lex} x^\eta$ for every term $x^\eta$ such that $x^\eta \mid x^\alpha $ and $x^{\alpha -\eta}
\notin J$.
 We will write $x^\alpha = x^\gamma \ast_J x^\delta$ to refer to this unique decomposition.
\item\label{it:borelProducts_ii} If  $x^\alpha \in J \setminus B_J$ and   $x_j = \min x^\alpha$, then $x^\alpha/x_j \in J$.
\item\label{it:borelProducts_iii} If $x^\beta\notin J$, while  $x^\delta x^\beta \in J$, then 
$x^\delta x^\beta = x^{\alpha}\ast_J x^{\delta'}$ with $x^\alpha \in B_J$ and $x^\delta >_{\mathtt{Lex}} x^{\delta'}$  (possibly 
$x^{\delta'}=1$). 
In particular, if $x_i x^\beta \in J$, 
then either $x_ix^\beta \in B_J$ or $x_i > \min x^\beta$.
\end{enumerate}
\end{proposition}

\begin{definition}\label{defF(s)}
 Let $J \in \mathbb{B}$  and $I$ be the ideal generated by a $J$-marked set $F$ in $A[\xx]$. We consider the following sets of polynomials:
\begin{itemize} 
\item$ F^{(s)} := \left\{ x^\delta f_\alpha\ \big\vert\ \deg \big(x^\delta f_\alpha\big) = t,\ f_\alpha \in F,\ \min x^\alpha \geq \max
x^\delta \ (\text{i.e.~} x^{\alpha+\delta}=x^\alpha*_Jx^\delta)\right\}$;
\item $\widehat{F}^{(s)}:=\left\{ x^\delta f_\alpha\ \big\vert\ \deg \big(x^\delta f_\alpha\big) = t,\ f_\alpha \in F,\ \min x^\alpha <
\max x^\delta \right\}$;
\item $\cN (J, I):= I \cap \langle \cN(J)\rangle$.  
\end{itemize} 
\end{definition}

Note that for $s=r$, we have  $F^{(r)}=F$, $\big\langle F^{(r)}\big\rangle=I_r$ and $\cN (J, I)_r=0$.

\begin{theorem}[{\cite[Theorems 1.7, 1.10]{LR2}}]\label{th:markedSetChar}
For $J \in \mathbb{B}$, let  $I$ be the ideal generated by  a $J$-marked set $F$ in $A[\xx]$.  Then, for every $s\geq r$,
\begin{enumerate}[(i)]
\item\label{it:markedSetChar_i} $I_s = \big\langle F^{(s)} \big\rangle + \big\langle \widehat{F}^{(s)} \big\rangle $;
\item\label{it:markedSetChar_ii} the $A$-module $\big\langle F^{(s)} \big\rangle$ is free of rank equal to $\vert F^{(s)} \vert =\rank( J_s)$;
\item\label{it:markedSetChar_iii} $I_s= \big\langle {F}^{(s)} \big\rangle \oplus  \cN (J, I)_s  $.
\end{enumerate}
Moreover,   the following conditions are equivalent:
\begin{enumerate}[(i)]\setcounter{enumi}{3}
\item\label{it:markedSetChar_iv} $F$ is a $J$-marked basis;
  \item\label{it:markedSetChar_v} for all $s\geq r$, $I_s=\big\langle {F}^{(s)}\big\rangle $;
\item\label{it:markedSetChar_vi}  $\cN(J,I)_{r+1}=0$;
\item\label{it:markedSetChar_vii}  $I_{r+1}=\big\langle {F}^{(r+1)}\big\rangle$;
\item\label{it:markedSetChar_viii} $\bigwedge^{Q+1} I_{r+1}=0$, where $Q:=\rk(J_{r+1})$.
\end{enumerate}
\end{theorem}
\begin{proof} This result is proved in a more general context in \cite{LR2}. We only  observe that   the conditions 
\lq\lq $\cN(J,I)_{s}=0$ and $I_{s}=\big\langle {F}^{(s)}\big\rangle $ 
for every
$s\leq \reg (J) +1$\rq\rq\ appearing in \cite{LR2} are equivalent  to \emph{(\ref{it:markedSetChar_vi})} and \emph{(\ref{it:markedSetChar_vii})}, since   
in the present hypotheses  $J$ is
generated in 
degree $r$ and $r$ is its regularity. 
With respect to \cite{LR2}, the only new item is \emph{(\ref{it:markedSetChar_viii})}, which  is obviously equivalent to 
\emph{(\ref{it:markedSetChar_vi})} and \emph{(\ref{it:markedSetChar_vii})}. 
In fact,  by  \emph{(\ref{it:markedSetChar_ii})} and  \emph{(\ref{it:markedSetChar_iii})} we have $I_{r+1}=\big\langle {F}^{(r+1)}\big\rangle \oplus 
\cN (J, I)_{r+1} $ 
and $\rank \big\langle {F}^{(r+1)}\big\rangle =\rk (J_{r+1})=Q$.
\end{proof}

\begin{corollary}\label{vuoto} Let $\mathcal I \in \mathcal E^{(0)}$ be such that $J(\mathcal I) \in \mathbb{B}$ and let $g\in \PGL$. Then:
\[ \underline{\mathbf{H}}_{\mathcal{I},g} \text{ is not empty } \Longleftrightarrow \ J(\mathcal I) \in   \mathbb{B}_{p(t)}.\]
Moreover, for $J=J(\mathcal I) \in \mathbb{B}_{p(t)}$ and any $k$-algebra $A$
\[
\underline{\mathbf{H}}_{\mathcal{I},g}(A)= \left\{  g\centerdot I \text{ s.t.~$I$ is generated by a } J\text{-marked basis in } A[\xx]  
\right\}.
\]
\end{corollary}
\begin{proof}  It is sufficient to prove the result for $g=\mathrm{id}$, i.e.~for  $\underline{\mathbf{H}}_{\mathcal{I}}$. 

Let $A$  be any $k$-algebra. If  $J=J(\mathcal I) \in  \mathbb{B}_{p(t)}$, then  $J\in \underline{\mathbf{G}}_{\mathcal{I}}(A)$ and the Hilbert
polynomial
of $\Proj(A[\xx]/J)$  is $p(t)$; hence $J \in \underline{\mathbf{H}}_{\mathcal{I}}(A)$. 

On the other hand, if $J=J(\mathcal I) \in \mathbb{B}\setminus  \mathbb{B}_{p(t)}$ and $I\in \underline{\mathbf{G}}_{\mathcal{I}}(A)$, then $I$ is generated 
by a $J$-marked set and $\rk (J_s) > q(s)$ for every  $s\gg 0$  (Lemma \ref{lemma:monomialsBorel}). By Theorem \ref{th:markedSetChar}, the
$A$-module $I_s$ 
contains a free submodule of rank equal to that of  $J_s$, hence  $I\notin \underline{\mathbf{H}}_{\mathcal{I}}(A)$. 

The second statement directly follows from Theorem \ref{th:markedSetChar}\emph{(\ref{it:markedSetChar_ii})} and the equivalence 
\emph{(\ref{it:markedSetChar_iv})}$\Leftrightarrow$\emph{(\ref{it:markedSetChar_v})}.
\end{proof}

In \cite{LR2} a functor $\underline{\mathbf{Mf}}_{J}\colon  \Rings  \rightarrow \Sets$ is defined for a strongly 
stable ideal $J$  in 
$\ZZ[x_0, \dots, x_n]$, by taking for a ring $A$
\[
\underline{\mathbf{Mf}}_{J}(A)= \left\{ \text{ideals } I \text{ generated by  } J\text{-marked bases in } A[x_0, \dots, x_n]  \right\}.
\]
Therefore, the open subfunctor $\underline{\mathbf{H}}_{\mathcal{I}}$ is the restriction  of $\underline{\mathbf{Mf}}_{J}$ to the sub-category 
$\Kalg$.

The marked functor $\underline{\mathbf{Mf}}_{J}$ is represented by a closed subscheme of the affine space $\AA^{D}_\ZZ$, for a suitable $D$. 
For the main features 
of $\underline{\mathbf{Mf}}_{J}$  and the proof of its representability see \cite{LR2}.  Here, we are only interested in the case $J\in \mathbb B$. 
Under this condition, 
any $J$-marked set has the shape \eqref{eq:markedset} and is uniquely determined by the list of $D=q(r)\cdot p(r)$ coefficients $c_{\alpha \gamma}$.  
Among the ideals generated by marked sets, 
those generated by marked bases are given for instance by  the closed condition \emph{(\ref{it:markedSetChar_viii})} (or  even \emph{(\ref{it:markedSetChar_vi})})
of Theorem \ref{th:markedSetChar}.

  \begin{theorem}\label{representable} The  Hilbert functor $\HilbFunctor{p(t)}{n}$ is the functor of points of a closed subscheme $\HilbScheme{p(t)}{n}$ of the 
  Grassmannian 
  $\GrassScheme{p}{N}$.
\end{theorem}
\begin{proof}
 By Proposition \ref{prop:openSubfunctors}, it  suffices to check the representability on an open cover of $\GrassFunctor{p}{N}$ and $\HilbFunctor{p(t)}{n}$:
 we choose 
 the Borel open cover
 (Definitions \ref{BorelcoverGrass} and \ref{BorelCoverHilb}). 
 For every $\mathcal I\in \mathcal E^{(0)}$ such that $J:=J(\mathcal I)\in \mathbb{B}_{p(t)}$ and  for every $g\in \PGL$,
 $\underline{\mathbf{H}}_{\mathcal{I},g}$ 
 is naturally isomorphic to $\underline{\mathbf{H}}_{\mathcal{I}}$. Moreover, $\underline{\mathbf{H}}_{\mathcal{I}}$  is the functor of
 points of  the $k$-scheme
  $\mathbf{H}_{\mathcal{I}}:={\mathbf{Mf}}_{J} \times_{\Spec \ZZ} \Spec k$.  
  Indeed, the scheme $\mathbf{H}_{\mathcal{I}}$ is the  subscheme of $\AA^D_k=\mathbf{G}_{\mathcal{I}}$ (where $D=p(r)\cdot q(r)$) defined by the closed equivalent 
  conditions    
of Theorem \ref{th:markedSetChar}. 
  Hence $\underline{\mathbf{H}}_{\mathcal{I}}$ is the functor of points of a closed subscheme $ {\mathbf{H}}_{\mathcal{I}}$ of $ {\mathbf{G}}_{\mathcal{I}}$.
  
  On the other hand if $J(\mathcal I) \in \mathbb{B} \setminus \mathbb{B}_{p(t)}$, then $\underline{\mathbf{H}}_{\mathcal{I}}$ is empty (Corollary \ref{vuoto}),
  hence it is the 
  functor of points of
  a closed subscheme  of   ${\mathbf{G}}_{\mathcal{I}}$. By Proposition \ref{prop:GrBorelSubFunctors} and the second part 
  of Proposition \ref{prop:openSubfunctors}, 
  we conclude that $\HilbFunctor{p(t)}{n}$ is the functor of points of a closed subscheme $\HilbScheme{p(t)}{n}$  of $\GrassScheme{p}{N}$.
\end{proof}
Next sessions are devoted to  describe how to determine  equations defining the Hilbert scheme $\HilbScheme{p(t)}{n}$ as subscheme of the 
Grassmannian $\GrassScheme{p}{N}$.


\section{Extensors and Pl\"ucker embedding}\label{sec:extensors}

In this section we consider any Grassmann functor, that we will denote by $\GrassFunctor{p}{N}$. In next sections, we will apply the tools developed to  the study of the  Hilbert functor and scheme. However, all the results of this section  hold true for every  $p$ and $N$, not only for
those obtained starting from an  Hilbert polynomial $p(t)$ of subschemes of $\PP^n$. 

In this section, we think at $\GrassFunctor{p}{N}(A)$  as presented  in  \eqref{eq:sottomodulo};  furthermore, our arguments allow us to restrict 
to  the open subfunctors $\underline{\mathbf{G}}_{\mathcal{I}}$, introduced in Section \ref{sec:Borel}. Thus,     the elements  of $ \GrassFunctor{p}{N}(A)$
we are mainly interested in are free submodules $L$  of $A^N$  of rank $q$, such that $A^N/L$ is free of rank $p$. 

We begin stating some well-know notions and results about exterior algebras.

\begin{definition}
Given a free $A$-module $M $, an \emph{extensor} of step $m$ in $M$ is an element of $\wedge^m M$ of the form 
$\mu_1 \wedge \cdots \wedge \mu_m$ with $\mu_1,\ldots,\mu_m$ in $M$.
\end{definition} 
Notice that $\mu_1 \wedge \cdots \wedge \mu_m$ vanishes whenever the submodule generated by $\mu_1,\ldots,\mu_m$ has rank lower than $m$.

\begin{lemma}\label{prop:flanders}
Let $\phi\colon P \rightarrow Q$ be a linear morphism of $A$-modules.
\begin{enumerate}[(i)]
\item\label{it:flanders_i} For any $m$, there exists a unique map $\wedge^m P \rightarrow \wedge^m Q$ such that
\[
p_1 \wedge \cdots \wedge p_m \mapsto \phi(p_1) \wedge \cdots \wedge \phi(p_m).
\]
We denote this morphism by $\phi^{(m)}$.
\item\label{it:flanders_ii} If $\phi$ is an isomorphism (resp.~surjective), then $\phi^{(m)}$ is an isomorphism 
(resp.~surjective) for every $m$.
\item\label{it:flanders_iii} If $\phi$ is injective and $P$ is free, then $\phi^{(m)}$ is injective
 for every $m$ \cite[Theorems 1, 8]{Flanders}. 
\item\label{it:flanders_iv} If $Q$ is   free   with basis $\{l_1,\ldots,l_s\}$, then for every $1\leqslant m \leqslant s$, the exterior 
algebra $\wedge^m Q$ is free of 
rank $\binom{s}{m}$ with basis $\{ l_{i_1}\wedge\cdots\wedge l_{i_m}\ \vert\ 1 \leqslant i_1 < \cdots < i_m \leqslant s \}$.
 
\noindent In particular, all the extensors of step $s=\rank Q$ associated to different bases of $Q$ are equal up to multiplication by an invertible element of $A$ \cite[Corollary A2.3]{eisen}.
\item\label{it:flanders_v}  If $M=P\oplus Q$, then  $\bigwedge^m(P\oplus Q)=\bigoplus_{r+s=m} \bigwedge^rP\otimes \bigwedge^s Q$.
\end{enumerate}
\end{lemma}

\begin{remark} As in the previous sections, $\a_1, \dots, \a_N$ is a fixed basis of the $A$-module $A^N$.  
We also  fix the isomorphism $\bigwedge^N A^N\simeq  A$ sending 
$\a_1 \wedge \cdots \wedge \a_N$ to $1_A$.  
For any $m$-multi-index $\mathcal J=(j_1, \dots, j_m)$,  we will denote by  $\a_{\mathcal J}$ the extensor
$ \a_{j_1}\wedge \dots \wedge \a_{j_m}$ of  $\wedge^m A^N$. By  Lemma \ref{prop:flanders}\emph{(\ref{it:flanders_iv})},
these extensors give a basis of $\wedge^m A^N$. We observe that  $\a_{\mathcal J}\wedge  \a_{\mathcal H}=0$ if 
$\mathcal H\cap {\mathcal J}\neq \emptyset $, while $\a_{\mathcal J}\wedge  \a_{\comp{\mathcal J}}
=\varepsilon_{{\mathcal J}} \a_1 \wedge \cdots \wedge \a_N$, where $\varepsilon_{{\mathcal J}}$ is  the signature of ${\mathcal J}{,}\comp{\mathcal J}$.  
Taking into account the fixed isomorphism,  we will simply write $\a_{\mathcal J}\wedge  \a_{\comp{\mathcal J}}
=\varepsilon_{{\mathcal J}}$.
\end{remark}

Every $A$-module  $L \in \underline{\mathbf{G}}_{\mathcal{I}}(A)$ has the    special free set of generators 
\begin{equation*}\mathcal B_{\mathcal I}(L):=\left\{ \b_{s}:=\a_{s} -\sum_{i\in \mathcal I} \gamma_{si}\, \a_{i}\ \bigg\vert \  s  \in \comp{\mathcal{I}} \right\}
\end{equation*}
described in Remark \ref{rk:precisazioni}. 
We will call it the ${\mathcal I}${\it{-marked set}} of $L$, extending the terminology we use in the special case of interest in this paper 
(Definition \ref{def:markedset})\footnote{The marked set $\mathcal B_{\mathcal I}(L)$ is in fact a basis for $L$; however, we do not call it
\lq\lq marked basis\rq\rq,  because in the case of a Grassmannian containing an Hilbert scheme, this terminology refers only to the points of the Hilbert scheme.}. 

\begin{definition}\label{Imarkedset}  
For every  $L \in \underline{\mathbf{G}}_{\mathcal{I}}(A)$  and  $\mathcal{S}=(s_1, \dots, s_m) \subset \comp{\mathcal{I}}$, 
we  denote by $\b_{\mathcal S}$ the extensor $\b_{s_1}\wedge  \dots \wedge  \b_{s_m}\in \wedge^m L$. 
The \emph{${\mathcal I}$-marked set} of $\wedge^m L$ is the  free set of generators  
\begin{equation*}\mathcal B_{\mathcal I}^{(m)}(L):=\left\{\b_{\mathcal S} \ \vert \ 
\mathcal S\subset \comp{\mathcal I}, \ \vert \mathcal S\vert =m \right\}. 
\end{equation*}
In particular, $\mathcal B_{\mathcal I}^{(1)}(L)=\mathcal B_{\mathcal I}(L)$.
\end{definition}

The aim  of the present section  is that of    determining a unified writing in terms of  the Pl\"ucker coordinates of  $ \GrassFunctor{p}{N}$ 
of a set  of generators of   $\wedge^m L$, where $1\leq m\leq q$ and $L\in \underline{\mathbf{G}}_{\mathcal{I}}(A)$. 
This set of generators will also contain  the   ${\mathcal I}$-marked set of  $\wedge^m L$. 

\medskip

   By  Lemma \ref{prop:flanders}\emph{(\ref{it:flanders_iii})}, there is a natural inclusion $  \wedge^mL \subseteq \wedge^m A^N$ for every
   $L \in \underline{\mathbf{G}}_{\mathcal{I}}(A)$. Hence,  every element $f\in \wedge^mL$ has a unique writing $f=\sum c_{\mathcal J} 
   \a_{\mathcal J}$, with coefficients $c_{\mathcal J}\in A$.

\begin{lemma}\label{lem:wedgeDec} Let $L \in \underline{\mathbf{G}}_{\mathcal{I}}(A)$.  
\begin{enumerate}[(i)]  
\item\label{it:wedgeDeci} If \ $\b_{\mathcal S} \in \mathcal B_{\mathcal I}^{(m)}(L)$ and   ${\mathcal K}:=\mathcal S \cup \mathcal I $, then
\begin{equation}\label{scritturabb}   \b_{\mathcal S}=
\a_{\mathcal S}+ \varepsilon^{\mathcal K}_{\mathcal S}\sum \varepsilon^{\mathcal K}_{\mathcal H}\,(\b_{\comp{\mathcal I}} \wedge \a_{ {\mathcal{K}}\setminus
{\mathcal H}} )\, \a_{\mathcal H } 
\end{equation}
 where  the sum is over    the  $m$-multi-indices $\mathcal H\neq   {\mathcal S}$ such that    $	 {\mathcal H}\subseteq   {\mathcal K}$, 
 and $\varepsilon^{\mathcal K}_{\mathcal H}$ is the signature of the permutation 
 $ \mathcal K\mapsto\mathcal H, {\mathcal K}\setminus {\mathcal H}$.
\item \label{it:wedgeDecii} If  $f=\sum c_{\mathcal J} \a_{\mathcal J}$ is  any  non-zero element  of  $  \wedge^mL\subset \wedge^mA^N$, then there is  
at least one non-zero  
coefficient $c_{\mathcal J}$ with $\mathcal J  \subset \comp{\mathcal I}$. 
\end{enumerate}
\end{lemma}
\begin{proof}
Up to a permutation,  we may assume that   $\comp{\mathcal K},{\mathcal S},\mathcal I=(1, \dots, N)$. 
Hence, $ \varepsilon^{{\mathcal K}}_{\mathcal S}=1$. 

\emph{(\ref{it:wedgeDeci})} We use the distributive law with $ \b_{s_j}=\a_{s_j} -\sum_{i\in \mathcal I} \gamma_{s_ji}\, \a_{i}$ and immediately see that
the coefficient of $\a_{\mathcal S}$ in $\b_{\mathcal S}$  is $1_A$, as  $\mathcal I\cap \mathcal S=\emptyset$, 
and the other extensors 
$\a_{\mathcal H}\neq \a_{\mathcal S}$ that   can appear with   non-zero coefficient  are those 
given in the statement.   As a consequence,   note that $ \b_{\mathcal S} \wedge \a_{\mathcal T}=0$ 
 if $	 {\mathcal T}$ is an $(N-m)$-multi-index  and $\comp{\mathcal T}\not\subseteq \mathcal K$, i.e.~$\mathcal T \not\supseteq   \comp{\mathcal K}  $.
 
Now we  prove the given formula for the coefficients, focusing  on each   $m$-multi-index 
$\mathcal H$.  Let us denote by $\gamma_{\mathcal  H} $ the  coefficient of $\a_{\mathcal H}$ in  $\b_{\mathcal S}$.
Applying again the distributive law on   $\a_{\comp{\mathcal K}}\wedge \b_{\mathcal S} \wedge \a_{\mathcal K \setminus  \mathcal H}$, 
the only non-zero summand is
$\gamma_{\mathcal  H }\,  (\a_{\comp{\mathcal K}}\wedge \a_{\mathcal H} \wedge \a_{\mathcal K \setminus  \mathcal H})=\gamma_{\mathcal  H }\, 
\varepsilon^{\mathcal K}_{\mathcal H}$, hence, 
$\gamma_{\mathcal  H }  =   \varepsilon^{\mathcal K}_{\mathcal H}\, (\a_{\comp{\mathcal K}}\wedge 
\b_{\mathcal S} \wedge \a_{\mathcal K \setminus  \mathcal H})$ and  
\begin{equation}\label{scritturab5}
\b_{\mathcal S}=\a_{\mathcal S}+ 
\sum\varepsilon^{\mathcal K}_{\mathcal H}\, (\a_{\comp{\mathcal K}}\wedge 
\b_{\mathcal S} \wedge \a_{\mathcal K \setminus  \mathcal H})\, \a_{\mathcal H } 
 \end{equation}
 with $\mathcal{H} \subseteq \mathcal{K}$, $ \vert\mathcal{H}\vert=m$, $\mathcal{H} \neq \mathcal{S}$. 

It remains to verify that 
$ (\a_{\comp{\mathcal K}}\wedge 
\b_{\mathcal S} \wedge \a_{\mathcal K \setminus  \mathcal H})=(\b_{\comp{\mathcal I}}\wedge \a_{\mathcal K \setminus  \mathcal H})$. 
  Applying \eqref{scritturab5}, we can write $\a_{\comp{\mathcal K}}=    \b_{\comp{\mathcal K}}-  \sum 
\gamma'_{\mathcal H'}\, \a_{\mathcal H' }$ where $\vert \mathcal H'\vert =\vert \comp{\mathcal K}\vert $ and  $\mathcal H'\neq   \comp{\mathcal K}$. 
We substitute and get
\[
 (\a_{\comp{\mathcal K}}\wedge 
\b_{\mathcal S} \wedge \a_{\mathcal K \setminus  \mathcal H})=(\b_{\comp{\mathcal K}}\wedge 
\b_{\mathcal S} \wedge \a_{\mathcal K \setminus  \mathcal H})
- \sum_{\mathcal H'} \gamma'_{\mathcal H' }\, (    \a_{\mathcal H'} \wedge \b_{{\mathcal S}} \wedge\a_{\mathcal{K}\setminus \mathcal H}).
\]
All the summands on $\mathcal H'$ vanish. Indeed, this is obvious if $\mathcal H'$ and   $ \mathcal{K}\setminus \mathcal H $ are not disjoint.
On the other hand,  if  they are and we denote by 
 $\mathcal T$ their union, then 
$\b_{{\mathcal S}} \wedge  \a_{\mathcal T}= 0$ since
 $\comp{\mathcal K}\cap {\mathcal{K}\setminus \mathcal H}=\emptyset $ and $\mathcal  H'\not\supseteq    \comp{\mathcal K}$.
Finally, 
$(\a_{\comp{\mathcal K}}\wedge 
\b_{\mathcal S} \wedge \a_{\mathcal K \setminus  \mathcal H})=(   \b_{\comp{\mathcal K}}\wedge \b_{{\mathcal S}} \wedge \a_{\mathcal{K}\setminus 
\mathcal H})= ( \b_{\comp{\mathcal I}}  \wedge 
\a_{\mathcal{K}\setminus \mathcal H}).$

\smallskip

\emph{(\ref{it:wedgeDecii})}  As $f \in \wedge^mL$, we can also write $f= \sum d_{\mathcal S}\,  \b_{\mathcal S}$, with  ${\mathcal S\subset 
\comp{\mathcal I}}$ and $d_{\mathcal S}\neq 0$.  By the previous item, $\a_{{\mathcal S}}$ appears only in the writing of $\b_{\mathcal S}$, hence 
its coefficient  $c_{\mathcal S}$  is $ d_{\mathcal S}\neq 0$.
\end{proof}

We would like to rewrite the coefficients appearing in the writing of the extensors $\b_{\mathcal S}$ given in (\ref{scritturabb}) 
in terms of  the Pl\"ucker  coordinates of $L$. Then, let us recall how they are defined.

\medskip

The projective space $\PP^E$ can be seen as the scheme representing the functor 
\[
\underline{\PP^E} \colon \Kalg  \rightarrow \Sets
\]
that associates to any
$k$-algebra $A$ the set
\begin{equation*}
\underline{\PP^E} (A) = \left\{\begin{array}{c}\text{isomorphism classes of epimorphisms } \\  \pi \colon A^{E+1} \rightarrow Q\text{ of 
locally free modules 
of rank }1\end{array} \right\}.
\end{equation*}
Hence, we can consider the natural transformation of functors $\underline{\mathscr{P}}\colon\GrassFunctor{p}{N} \rightarrow
 \underline{\PP^{E}}$ given by:
\begin{equation*}
\underline{\mathscr{P}}(A)\colon \ \left(\phi_L\colon  A^N \xrightarrow{\pi} A^N/L\right) \in \GrassFunctor{p}{N}(A) 
\quad \longmapsto \quad \left(\phi_L^{(p)} \colon \wedge^p A^N \xrightarrow{\pi^{(p)}} 
\wedge^p ( A^N/L ) \right) \in 
\underline{\PP^{E}}(A)
\end{equation*}
where   $\wedge^p A^N $ is free of rank $ \binom{N}{p}=E+1 $ and $\wedge^p A^N/L $ is locally free of rank 1. 

The collection of open subfunctors
$\underline{\mathbf{G}}_{\mathcal{I}}$ of Proposition \ref{prop:openCoverGrass} is exactly that  induced by the transformation
$\underline{\mathscr{P}}$
and the standard affine cover of the projective space $\PP^E$
corresponding to the  basis   $\{\a_{\mathcal J} \ \vert \ \mathcal J \in \mathcal E^{(0)}\}$ 
of $\wedge^p A^N$. 

We denote by $\mathbf{\Delta}$ the variables of $\PP^E$ and we index them  using the  multi-indices $\mathcal I \in \mathcal{E}^{(0)}$ 
so that
 ${\mathbf{G}}_{\mathcal{I}} $ be the open subscheme of the Grassmannian 
defined by the condition $\mathbf{\Delta}_\mathcal{I} \neq 0$.
The Grassmannian $\GrassScheme{p}{N}$ is a closed subscheme of $\PP^E = \Proj k[\mathbf{\Delta}]$ 
defined by the Pl\"ucker relations, that are generated  
by homogeneous polynomials of degree $2$: we will denote by $k[\DD]$   the coordinate ring of $\GrassScheme{p}{N}$,   i.e.~the quotient of $k[\mathbf{\Delta}]$ 
under the Pl\"ucker relations, so that $\GrassScheme{p}{N} = \Proj k[\DD] \subset \PP^E = \Proj k[\mathbf{\Delta}]$ (see for instance \cite{KleimanLaksov}). 

We can also associate Pl\"ucker coordinates to each module $L\in \underline{\mathbf{G}}_{\mathcal{I}}(A)$. 
Upon 
fixing an isomorphism $i\colon \wedge^p (A^N/L) \simeq  A$, $\underline{\mathscr{P}}(A)(L)$ can be seen as the map $i\circ\phi_L^{(p)}$  or, equivalently,
as the function
\[
i\circ \phi_L^{(p)}\colon  \left\{\a_{\mathcal J} \ \big\vert \ \mathcal J \in \mathcal E^{(0)}\right\}\rightarrow A \ \text{ given by  }  \  
\a_{\mathcal J} \mapsto  {\Delta}_{\mathcal{J}}(L):=i\big(\phi_L^{(p)}(\a_{\mathcal J})). \quad
\]

Since two  isomorphisms
$i,i'\colon \wedge^p A^N/L \rightarrow A$ only differ by the multiplication by a unit $u\in A$,  the  Pl\"ucker coordinates of $L$ are 
defined up to invertible
elements in $A$.  

  By definition of  $\underline{\mathbf{G}}_{\mathcal{I}}(A) $, we have  the decomposition as direct sum 
  $A^N=L\oplus \langle \a_i \ \vert\ i\in \mathcal I \rangle$, so that $ \phi_L^{(p)}$ factors through  $\wedge^pA^N \rightarrow \wedge^N  A^N\rightarrow \wedge^p (A^N/L)$ given by $\a_{\mathcal J}\mapsto \b_{\comp{\mathcal I}}\wedge \a_{\mathcal J}\mapsto \overline{\a_{\mathcal J}}$, where 
 $\b_{\comp{\mathcal I}}$   is  the only element 
 of the $\mathcal I$-marked set $ \mathcal B_{\mathcal I}^{(q)}(L)$ of $\wedge^qL$.

Hence, the Pl\"ucker coordinates of $L$ are
\begin{equation}\label{lem:wedgeDecPlucker}
\left( {\Delta}_{\mathcal{J}}(L)= \b_{\comp{\mathcal I}} \wedge \a_{\mathcal J}\  \big\vert \  \mathcal J\in \mathcal E^{(0)} \right).
\end{equation}

We\hfill identify\hfill $\b_{\comp{\mathcal I}} \wedge \a_{\mathcal J}$\hfill  with\hfill elements\hfill of\hfill $A$\hfill by\hfill fixing\hfill the\hfill isomorphisms\hfill $\wedge^N A^N\simeq A$\hfill and\\
$i\colon \wedge^p A^N/L \rightarrow A$.  For the first one we fixed that sending $\a_{(1, \dots, N)} $ to 1; if we choose   
$i\colon  \overline{\a_{\mathcal I}}  \mapsto 1$, then   \eqref{lem:wedgeDecPlucker} gives the representative of the  Pl\"ucker coordinates with
$ {\Delta}_{\mathcal{I}}(L)=1$. Indeed,   in our setting $\varepsilon_{\comp{\mathcal I}}=1$ and
$\b_{\comp{\mathcal I}} \wedge \a_{\mathcal I}=\a_{\comp{\mathcal I}} \wedge \a_{\mathcal I }=\a_{(1, \dots, N)}$ by Lemma \ref{lem:wedgeDec}.

 Therefore, Pl\"ucker coordinates of $L$  can  be obtained as the maximal minors of the $q\times N$ 
 matrix whose rows contain 
the elements of $\mathcal B_{\mathcal I}(L)$. 
 More precisely,   $ \Delta_{\mathcal J}(L)$ is  the   minor
corresponding to the columns with indices in $\comp{\mathcal{J}}$,  up to  a sign   given by  the signature  $\varepsilon_{\mathcal{J}}$.

\medskip

Using \eqref{lem:wedgeDecPlucker} we can finally  rewrite the coefficients appearing in (\ref{scritturabb}) in terms 
of the  Pl\"ucker coordinates of $L$.

\begin{corollary}\label{prop:genDelta}
Let $L \in \underline{\mathbf{G}}_{\mathcal{I}}(A)$, $\b_{\mathcal S}$ be  any extensor in
$\mathcal B^{(m)}_{\mathcal I}(L)$  and   $
 {\mathcal K}:=\mathcal S \cup \mathcal I $.
Then  
\begin{equation*}\label{scritturab2} \varepsilon^{\mathcal K}_{\mathcal S}\, \Delta_{ {\mathcal{I}}}(L)\, \b_{\mathcal S}= 
 \varepsilon^{\mathcal K}_{\mathcal S}\, \Delta_{ {\mathcal{I}}}(L)\, \a_{\mathcal S}+ 
\sum 
\  \varepsilon^{\mathcal K}_{\mathcal H}\,  
\Delta_{{\mathcal{K}}\setminus \mathcal H}(L)\, \a_{\mathcal H } 
\end{equation*}
where the sum is over  the $m$-multi-indices $\mathcal H$   such that $\mathcal H\subseteq 
 {\mathcal K}$, $\mathcal H\neq \mathcal S$.
\end{corollary}

\begin{definition}\label{defB}
For every $1\leq m\leq q$, we define the following subset of  $ k[\DD]^N$:
\begin{equation*}
\mathcal{B}^{(m)}:=
\left\{\delta^{(m)}_{\mathcal{K}}:= \sum_{\begin{subarray}{c} \mathcal{H} \subseteq \mathcal{K} \\ \vert\mathcal{H}\vert=m \end{subarray}}
\varepsilon^{\mathcal K}_{\mathcal H}\, {\Delta}_{\mathcal{K}\setminus\mathcal{H}}\, \a_{\mathcal H}\ \Bigg\vert\ \mathcal{K} 
\in \mathcal{E}^{(m)} \right\}.
\end{equation*}
Moreover,  for every $\mathcal{I} \in \mathcal{E}^{(0)}$, we define $  \mathcal{B}_{\mathcal{I}}^{(m)}:=
\left\{\delta^{(m)}_{\mathcal{K}}\ \big\vert\ \mathcal{K} 
\in \mathcal{E}^{(m)},\ \mathcal{K} \supseteq \mathcal I \right\}$.
\end{definition}

\begin{remark} \label{sommadirettaDelta}
 For every $m$-multi-index $\mathcal S$ that does not intersect $\mathcal I$, $\a_{\mathcal S}$ appears in a single element of  $  \mathcal{B}_{\mathcal{I}}^{(m)}$, 
  the one with index $\mathcal K=\mathcal I\cup \mathcal S$. Moreover, 
  $\delta^{(m)}_{\mathcal{K}}-\varepsilon^{\mathcal K}_{\mathcal S}\,{\Delta}_{\mathcal{I}}\, \a_{\mathcal S}$ is the 
sum 
 $\sum \varepsilon^{\mathcal K}_{\mathcal H}\, {\Delta}_{\mathcal{K}\setminus\mathcal{H}}\, \a_{\mathcal H}$, 
 where all the  $m$-multi-indices $\mathcal H$   intersect $\mathcal I$. 
 
Hence, for every element  $ f\in \wedge^m k[\DD]^N $ we can write $ {\Delta}_{\mathcal{I}}f$ as a sum $f_1+f_2$ with 
 $f_1\in \langle \mathcal B^{(m)}_{\mathcal I} \rangle$ and $f_2\in \langle 
 \a_{\mathcal H} \text{ s.t.~}\vert\mathcal H\vert=m \text{ and }  \mathcal H\cap \mathcal I\neq \emptyset \rangle$.  
\end{remark}

We will now evaluate the elements   $\delta^{(m)}_{\mathcal{K}}$ at  $L \in \underline{\mathbf{G}}_{\mathcal{I}}(A)$; of course, such evaluations are defined
only up 
to   units of $A$. Through evaluation at $L$, we can also see that   the  notations of Definition \ref{defB} are consistent with those introduced 
in Definition \ref{Imarkedset} 
to denote $\mathcal I$-marked  sets of $\wedge^m L$. 

\begin{theorem}\label{cor:costbase} Let $\mathcal I\in \mathcal E^{(0)} $, $A$ be a  $k$-algebra and  $L $ be a module in 
$ \underline{\mathbf{G}}_{\mathcal{I}}(A)$. Then,  for every $1\leq m\leq q$,
\begin{enumerate}[(i)]
\item\label{it:costbase_i}  the evaluation of    $\mathcal B^{(m)}_{\mathcal I}$ at $L$
is the $\mathcal I$-marked set $\mathcal B^{(m)}_{\mathcal I}(L)$
 of $\wedge^m L$;
 \item\label{it:costbase_ii} the evaluation $\mathcal B^{(m)}(L)$ of  $\mathcal B^{(m)}$ at $L $ is a set of generators of $\wedge^m L$.
\end{enumerate}
\end{theorem}

\begin{proof}
\emph{(\ref{it:costbase_i})} Let $\mathcal K$ be  a $(q+m)$-multi-index containing $\mathcal I$ and let $\mathcal S=\mathcal K\setminus \mathcal I$.
As a straightforward consequence of Lemma \ref{lem:wedgeDec} and Corollary \ref{prop:genDelta}, we see that 
 $ \delta^{(m)}_{\mathcal{K}}(L)$ is equal (up to units of $A$) to the element  $\b_{\mathcal S}$ of the $\mathcal I$-marked set of $\wedge^m L$.
 Note that for $L \in \underline{\mathbf{G}}_{\mathcal{I}}(A)$ and 
 $\mathcal H=\mathcal S$ we 
 may set ${\Delta}_{\mathcal{K}\setminus\mathcal{S}}(L)={\Delta}_{\mathcal{I}}(L)=1$.
 
 \smallskip
 
\emph{(\ref{it:costbase_ii})} By the previous item, it  suffices to prove that $\mathcal B^{(m)}(L)\subset \wedge^m L$.

 Let us consider any 
 $ \delta^{(m)}_{\mathcal{K'}}\in \mathcal B^{(m)} $ and    write 
 $ {\Delta}_{\mathcal{I}}\, \delta^{(m)}_{\mathcal{K'}}= \delta_1+\delta_2$ 
  as in Remark \ref{sommadirettaDelta} with 
 $\delta_1\in \langle \mathcal B^{(m)}_{\mathcal I} \rangle$ and $\delta_2\in \langle 
 \a_{\mathcal H} \text{ s.t.~} \vert \mathcal H\vert=m \text{ and } \mathcal H\cap \mathcal I\neq \emptyset \rangle$.   
 Under our assumption,  ${\Delta}_{\mathcal{I}}(L)$ is a unit in $A$; therefore, we need to prove that $\delta_2(L)=0$.

If there is a $p$-multi-index $\mathcal I' \subset \mathcal K'$ such that  $L \in \underline{\mathbf{G}}_{\mathcal{I'}}(A)$, then it follows by
\emph{(\ref{it:costbase_i})} that
$\delta^{(m)}_{\mathcal{K'}}(L)\in \mathcal B^{(m)}_{\mathcal I'}(L)\subset \wedge^m L$, so that also 
$\delta_2(L) \in \wedge^m L$ and we get 
$\delta_2(L)=0$ by Lemma \ref{lem:wedgeDec}\emph{(\ref{it:wedgeDecii})}.

Therefore, $\delta_2$ vanishes over the non-empty open subfunctor  $ \underline{\mathbf{G}}_{\mathcal{I}}(A)\cap  \underline{\mathbf{G}}_{\mathcal{I'}}(A)$
of the Grassmann functor, hence it vanish on $\GrassScheme{p}{N}$.
\end{proof}

We will use the results of this section  in order to compute equations defining \emph{globally} the Hilbert scheme  as subscheme of the Grassmannian,
starting from those defining   ${\mathbf{H}}_{\mathcal{I}}$ in ${\mathbf{G}}_{\mathcal{I}}$. 
Then in the following the elements of basis $\a_1, \dots, \a_N$ of $A^N$  will correspond to the terms  $x^{\alpha(1)},
\dots,x^{\alpha(N)}$ in $k[\xx]_r$. We can reformulate Theorem \ref{cor:costbase} in this special setting.

\begin{theorem}\label{universal} The universal family $\mathcal{F} \hookrightarrow \PP^n \times \GrassScheme{p}{N} \rightarrow \GrassScheme{p}{N}$ 
parameterized by the 
Grassmannian, given in \eqref{eq:universalfamily}, is generated by the set of bi-homogeneous elements  in $k[\DD , \xx]$
\begin{equation*} 
\left\{  \delta^{(1)}_{\mathcal{K}} = \sum_{ h\in\mathcal{K} } \varepsilon^{\, \mathcal K}_{\{h\}}\, \Delta_{\mathcal{K}\setminus\{h\}} \,
x^{\alpha(h)} \quad \Big\vert \quad \forall\ \mathcal{K} \in \mathcal{E}^{(1)} \right\}
\end{equation*}
and the $m$-th exterior power of the universal element is generated by
\begin{equation*}
\left\{  \delta^{(m)}_{\mathcal{K}} = \sum_{\begin{subarray}{c}\mathcal{H}\subset\mathcal{K} \\ \vert\mathcal{H}\vert = 
m \end{subarray} } \varepsilon^{ \mathcal K}_{\mathcal{H}}\,\Delta_{\mathcal{K}\setminus \mathcal{H}} \,
x^{\alpha(h_1)} \wedge \cdots \wedge x^{\alpha(h_m)} \quad \Bigg\vert \quad \forall\ \mathcal{K} \in \mathcal{E}^{(m)}\right\}.
\end{equation*}
\end{theorem}

\section{Equations}\label{sec:BLMR}
 In this section we will obtain global equations defining the Hilbert scheme. In particular, the new set of equations has degree lower than the other known equations. Towards this aim we need to refine some results of Section \ref{sec:representability}, in particular  Theorem \ref{th:markedSetChar}.

 These results concern any  $J$-marked set $F$, where $J$ is a Borel-fixed ideal generated by $q$ terms of degree $r$; we do not assume that 
 the Hilbert polynomial $p_J(t)$ of $A[\xx]/J$ is $p(t)$.  However,   we know  that   $r$ is the regularity of $J$ and, by Lemma \ref{lemma:monomialsBorel}, 
  $k[x_{d+1},\ldots, x_{n}]_{\geq r} \subset J$ and  $\deg p_J(t)\leq  d=\deg p(t)$; hence 
$\cN(J)_{\geq r} \subset (x_0, \dots, x_d)$.  In particular,   the support of every polynomial 
 $f_\alpha\in F$ is contained in  $(x_0, \dots, x_d)$, except for only one possible term, the head term $\Ht(f_\alpha)=x^\alpha$.

\begin{definition}
Let $J\in \mathbb B$ and let $I \subset A[\xx]$ be an ideal generated by a $J$-marked set $F$. Making 
reference (and in addition) to Definition \ref{defF(s)},   we set:
\begin{itemize}
\item $F':=\left\{x_i f_\alpha \in F^{(r+1)} \ \vert \  i=d+1,\dots, n\right\}=F^{(r+1)} \setminus (x_0, \dots, x_d)$;   
\item $F'':=\left\{x_i f_\alpha \in F^{(r+1)}  \ \vert \   i=0, \dots, d\right\}= F^{(r+1)} \cap(x_0, \dots, x_d)$;
\item $S:=\left\{ x_j f_\beta-x_if_\alpha \ \vert \ \forall  x_j f_\beta \in \widehat{F}^{(r+1)},   x_i f_\alpha \in F^{(r+1)} \text{ s.t.~} 
x_jx^\beta=x_ix^\alpha \right\}$;
 \item  $ q':=\dim_k k[x_{d+1},\ldots, x_{n}]_{r+1}$;  
\item $   q'':= q(r+1) -q'$;
\item  $I'':= I_{r+1} \cap (x_0, \dots, x_d)$;
\item $I^{(1)}:=\langle x_hI_r \ \vert \ \forall\ h=0, \dots, d  \rangle \subseteq I''$.
\end{itemize} 
\end{definition}

\begin{theorem}\label{thm:rifatto} Let $J\in \mathbb B$ and $I \subset A[\xx]$ be an ideal generated by a $J$-marked set $F$. 
Then, 
\begin{enumerate}[(i)]
 \item\label{it:rifatto_i}   $\langle   {F}'\rangle $ is a  free $A$-module of  rank  $q'$; 
 \item\label{it:rifatto_ii}  $\langle   {F}''\rangle $ is a  free $A$-module contained in $ I^{(1)}$ of  rank  $\geq q''$;
 \item\label{sommadiretta_2}   $I_{r+1}=\langle {F'}\rangle \oplus I''$;
 \item\label{sommadiretta_3}   $I''= \langle {F''} \rangle \oplus   \cN(J,I)_{r+1}= \langle {F''} \rangle + \langle S \rangle $.
\end{enumerate}
Moreover, the following conditions are equivalent:
\begin{enumerate}[(i)]\setcounter{enumi}{4}
\item\label{it:rifatto_v}  $J\in\mathbb B_{p(t)}$ and  $F$ is a $J$-marked basis;
\item \label{it:rifatto_vi} $\wedge^{q(r+1)+1} I_{r+1}=0$;
\item \label{it:rifatto_vii} $\wedge^{q''+1} I''=0$;
\item \label{it:rifatto_viii}  $\wedge^{q''+1}I^{(1)}=0$ and $\big(\wedge^{q''} I^{(1)}\big) \wedge I''=0$.
\end{enumerate}
\end{theorem}

\begin{proof} 
\emph{(\ref{it:rifatto_i})}  It is sufficient to recall that   $F' $ is a subset of the set of linearly independent polynomials $F^{(r+1)}$, 
hence the $A$-module $\langle   
{F}'\rangle $ is free of rank equal to $\vert F'\vert $. Moreover  $\vert F'\vert=q'$ by Lemma \ref{lemma:monomialsBorel}.

\smallskip

\emph{(\ref{it:rifatto_ii})} We can prove  that  $\langle   {F}''\rangle $ is free with rank $\vert F''\vert$ by the same argument used for
\emph{(\ref{it:rifatto_i})}.
Moreover, by definition and Lemma \ref{lemma:monomialsBorel},  $\langle{F}''\rangle =\vert F^{(r+1)}\vert -\vert F'\vert= \rk(J_{r+1})- \vert F'\vert \geq q(r+1)-q'$.

\smallskip

\emph{(\ref{sommadiretta_2})},\emph{(\ref{sommadiretta_3})}   We obtain  the equality  
$I_{r+1}= \langle {F'}\rangle \oplus \langle {F''} \rangle \oplus   \cN(J,I)_{r+1} $ as a  consequence of
 Theorem \ref{th:markedSetChar}\emph{(\ref{it:markedSetChar_iii})} and the fact that $F^{(r+1)}$ is the disjoint union of $F'$ and $F''$. 
 It is obvious by the definition that  
$ \langle {F''} \rangle \oplus  \cN(J,I)_{r+1} \subseteq I''$.
Then,  to prove      $I_{r+1}=\langle {F'}\rangle \oplus I''$ 
it   suffices to verify that  the sum $\langle {F'}\rangle + I''$ is direct. If $h$ is any element $h=\sum d_{i\alpha}x_if_\alpha \in 
\langle {F'}\rangle$ with  $d_{i\alpha}\in A$, $d_{i\alpha}\neq 0$,  then $x_ix^\alpha \in \supp(h)$, since the head terms of the monic marked polynomials  
$x_if_\alpha\in F'$ are distinct terms in $k[x_{d+1}, \dots,x_n]_{r+1}$, while   $x_if_\alpha-x_ix^\alpha \in (x_0, \dots, x_n)$.
 Therefore, we also get $I''=   \langle {F''} \rangle \oplus  \cN(J,I)_{r+1}$,

Let us consider the set of generators $F'\cup F''\cup \widehat{F}^{(r+1)}$  of the $A$-module $I_{r+1}$. For every element  
$x_j f_\beta\in \widehat{F}^{(r+1)}$,  we can find an element 
 $x_if_{\alpha}\in F'\cup F''$ such that  $x_ix^{\alpha}=x_jx^\beta$ and  $h_{j\beta}:=x_jf_\beta -x_if_\alpha\in S$. Then, we get 
 a new set of generators replacing $\widehat{F}^{(r+1)}$ by $S$.  The union of the three sets $F'$, $F''$ and $S$ generates the  
 $A$-module $I_{r+1}$ and, in particular,  $F''\cup S$ generates $I''$, since $S\subseteq      I''$.  

\smallskip

\emph{(\ref{it:rifatto_v})}$\Leftrightarrow$\emph{(\ref{it:rifatto_vi})} If $J\in  \mathbb{B}_{p(t)}$, then the statement is given 
by Theorem \ref{th:markedSetChar}\emph{(\ref{it:markedSetChar_iv})}$\Leftrightarrow$\emph{(\ref{it:markedSetChar_viii})}, as $\rk(J_{r+1})=q(r+1)$.
On the other hand, if $J\notin  \mathbb{B}_{p(t)}$, then by Gotzmann's Persistence Theorem we have $\rk(J_{r+1})>   q(r+1)$, so that  $\wedge^{q(r+1)+1}
I_{r+1}\neq 0$ by Theorem \ref{th:markedSetChar}\emph{(\ref{it:markedSetChar_ii})}. 

\smallskip

\emph{(\ref{it:rifatto_vi})}$\Leftrightarrow$\emph{(\ref{it:rifatto_vii})}$\Leftrightarrow$\emph{(\ref{it:rifatto_viii})} are 
straightforward consequences of previous items.
\end{proof}

\begin{proposition}\label{solo2} In the setting of Theorem  \ref{thm:rifatto},   let $ B$ be any  set of polynomials of $I_r$ 
containing $F$ and consider the following two subsets of $I_{r+1}$:  
\begin{enumerate}[1)]
\item\label{it:1}   $\bigcup_{i=0}^s x_i B$;
\item\label{it:2}   $\{ x_i f -x_{j}g \ \vert \ \forall\  f, g \in  B  \text{ such that  } x_i f-x_jg \in (x_0, \dots, x_d)\}$.
\end{enumerate}
For $s=n$  the elements in \ref{it:1}) generate $I_{r+1}$,  while for  $s=d$  they generate  $I^{(1)}$.  Moreover, the first set 
for $s=d$ and the second set generate $I''$. 
\end{proposition}
\begin{proof}   The first and second assertions are straightforward by the definitions    of $I_{r+1}$ and  $I^{(1)}$.
 For the latter one, we observe that the polynomials in these two  sets are contained in $I''=(F)_{r+1}\cap(x_0, \dots, x_d)$.  Thus,  it suffices 
 to prove the statement  in the case $B=F$.

By Theorem \ref{thm:rifatto}\emph{(\ref{sommadiretta_3})},  the $A$-module $I''$ is generated by  $F'' \cup S$. Obviously, $F''$ is contained 
in the set given in \emph{\ref{it:1})}. Moreover,  $S$ in contained in the set  given in \emph{\ref{it:2})}. Indeed, 
by definition of $J$-marked set and Lemma \ref{lemma:monomialsBorel},  for every  $f_\alpha \in F$ we have  
$f_\alpha -x^\alpha \subset \langle \cN(J)\rangle_r\subset (x_0, \dots, x_d)$. 
Then,  $f_\alpha\in (x_0, \dots, x_d)$ if, and only if, 
 $x^\alpha \in (x_0,\dots, x_d)$. 
\end{proof}

\begin{remark}\label{rk:generat} For every ideal $I\in  \underline{\mathbf{G}}_{\mathcal{I}}(A)$ with $J(\mathcal I)\in \mathbb B$, we will apply 
the previous  results considering $J(\mathcal I)$ as $J$  and   the set of 
generators  $\mathcal{B}^{(1)}(I)$ (where $I$ stands for $I_r $) as $B$. 
Note that  $\mathcal{B}^{(1)}(I)$   contains the    $\mathcal I$-marked set  $\mathcal{B}^{(1)}_{\mathcal I}(I)$, which is monic
since $\Delta_{\mathcal I}(I)$ is a unit in $A$ and we may set   
$\Delta_{\mathcal I}(I)=1$.

In order to apply to $I$ the equivalent conditions \emph{(\ref{it:rifatto_v})},\dots,\emph{(\ref{it:rifatto_viii})} of 
Theorem \ref{thm:rifatto} we need to consider exterior products of the type
$ \wedge^m \langle x_0 I_r, \dots , x_sI_r\rangle $ for some integers $1\leq m\leq q(r+1)+1$ and $0\leq s\leq  N$. 
The  set of generators for this module we  use is 
\begin{equation*}
\left\{  \bigwedge_{\begin{subarray}{c} 0 \leq i \leq s \\ m_i > 0\end{subarray}} x_i
\delta^{(m_i)}_{\mathcal{K}_i}(I)
\ \ \Bigg\vert\  \ \forall\	\delta^{(m_i)}_{\mathcal{K}_i} \in \mathcal{B}^{(m_i)} \text{ s.t.~}  \sum m_i=m  \right\}.
\end{equation*}
This set is obtained considering the decomposition of $\wedge^m \langle x_0 I_r, \dots , x_sI_r\rangle$ as the sum of the submodules
$(x_0\wedge^{m_0}I_r) \wedge \dots \wedge (x_s \wedge^{m_s}I_r) $ 
over the sequences of non-negative  integers $(m_0, \dots, m_s)$ with sum $m$.
Note that in this   writing we assume that the $i$-th piece $x_i\wedge^{m_i}I_r$ is missing  whenever $m_i=0$;
 the number of factors  is at most $s$
and the maximum is reached only if all the integers $m_i$ are positive.
\end{remark}

\medskip
We are now able to exhibit the ideal $ \mathfrak{H}$ in    the ring of Pl\"ucker coordinates $k[\Delta]$ 
that globally   defines the Hilbert scheme as a subscheme of the Grassmannian. First, we set 
\begin{eqnarray}
\label{eq:part1}&&
\mathfrak{h}_{1}:=\xcoeff\left\{  \bigwedge_{\begin{subarray}{c} 0 \leq i \leq d \\ m_i > 0\end{subarray}} x_i
\delta^{(m_i)}_{\mathcal{K}_i}  
\ \bigg\vert  \ \forall\	\delta^{(m_i)}_{\mathcal{K}_i} \in \mathcal{B}^{(m_i)}  \text{ s.t.~}    \sum m_i=q''+1   \right\}\\
 &&\label{eq:part2}
 \mathfrak{h}_{2}:=\xcoeff\left\{   \left(\bigwedge_{\begin{subarray}{c} 0 \leq i \leq d \\ m_i > 0\end{subarray}} 
 x_i \delta^{(m_i)}_{\mathcal{K}_i}\right) \wedge \left(x_j\delta^{(1)}_{\mathcal{H}} \pm  x_{k}\delta^{(1)}_{\overline{\mathcal{H}}}\right) 
 \ \Bigg\vert
   \begin{array}{l} \forall\	\delta^{(m_i)}_{\mathcal{K}_i} \in \mathcal{B}^{(m_i)}  \text{ s.t.~}   \sum m_i=q''\\
   \forall\ x_j\delta^{(1)}_{\mathcal{H}} \pm x_{k}\delta^{(1)}_{\overline{\mathcal{H}}} \in \mathcal{W}
  \end{array} 
  \right\}\qquad
\end{eqnarray}
where 
$\mathcal W$ is the set of polynomials $x_j\delta^{(1)}_{\mathcal{H}} \pm x_{k}\delta^{(1)}_{\overline{\mathcal{H}}}$ such that 
\begin{itemize}
\item $\mathcal{H} = (\mathcal{H}\cap\overline{\mathcal{H}}) \cup \{h\}$ and $\overline{\mathcal{H}} = (\mathcal{H}\cap\overline{\mathcal{H}}) \cup \{\overline{h}\}$, i.e.~the polynomial $\delta^{(1)}_{\mathcal{H}}$ contains the term $\Delta_{\mathcal{H}\cap\overline{\mathcal{H}}}x^{\alpha(h)}$ and $\delta^{(1)}_{\overline{\mathcal{H}}}$ contains $\Delta_{\mathcal{H}\cap\overline{\mathcal{H}}}x^{\alpha(\overline h)}$;
\item the pair $(x_j,-x_k)$ is a syzygy for the monomials $x^{\alpha(h)}$ and $x^{\alpha(\overline{h})}$, i.e. $x_j x^{\alpha(h)} = x_k x^{\alpha(\overline{h})}$ and the sign $\pm$ is chosen in order for the terms $\Delta_{\mathcal{H}\cap\overline{\mathcal{H}}} x_j x^{\alpha(h)}$ and $\Delta_{\mathcal{H}\cap\overline{\mathcal{H}}} x_k x^{\alpha(\overline{h})}$ to cancel;
\item $\xx \text{-supp}(x_j\delta^{(1)}_{\mathcal{H}} \pm x_{k}\delta^{(1)}_{\overline{\mathcal{H}}})\subset (x_0, \dots, x_d)$.
\end{itemize}

\medskip

Moreover, we set  $\mathfrak{h}:=\mathfrak{h}_{1} \cup \mathfrak{h}_{2} $ and consider for every $g\in \PGL$ the set of equations
$g\centerdot \mathfrak{h}$ 
obtained  by   the action of $g$  on the elements of  $\mathfrak{h}$.  Finally we define the ideal
  \begin{equation*}
  \mathfrak{H}:=\left(\bigcup_{g \in \PGL} g\centerdot \mathfrak{h}\right). 
  \end{equation*}

\begin{theorem}\label{th:mainTheorem} Let $p(t)$ be an admissible Hilbert polynomial for subschemes of $\PP^n$ of degree $d$. 

The homogeneous ideal $\mathfrak H$  in the ring of Pl\"ucker coordinates  $k[\DD]$ of the Pl\"ucker 
embedding $\GrassScheme{p}{N} \hookrightarrow \PP^E$  is generated in  degree 
 $\leq d+2$  and  defines  
 $\HilbScheme{p(t)}{n}$ 
 as a closed subscheme of $\GrassScheme{p}{N}$.
\end{theorem}

\begin{proof} 
By  definition,  $\mathfrak H$  is the smallest ideal in $k[\DD]$ that  contains  the union of the  two sets
of equations $\mathfrak{h}_{1}$ and $\mathfrak{h}_{2}$, given in 
 \eqref{eq:part1} and \eqref{eq:part2}, 
  and    is invariant by the action of $ \PGL$. Since the action of $ \PGL$ does not
 modify the degree of polynomials, in order to prove  the first part of the statement it  suffices to recall that 
 each $\delta^{(m)}_{\mathcal{K}}$ is linear in the Pl\"ucker coordinates (Theorem \ref{universal}); hence,   
the degree of each polynomial 
in \eqref{eq:part1} is at most $d+1$ and the degree of each polynomial in 
\eqref{eq:part2} is  at most $d+2$. In both cases equality is achieved only when all the integers $m_i$ are strictly positive. 

\medskip

For convenience,  we  denote by $\mathcal{Z} $ the subscheme of
$\GrassScheme{p}{N}$ defined by $\mathfrak{H} $ and by  $\mathfrak{D}$ the saturated ideal in $k[\Delta]$ that defines $\HilbScheme{p(t) }{n}$
 as a closed subscheme of $\GrassScheme{p}{N}$.  We  have to prove that $\mathcal{Z} =\HilbScheme{p(t)}{n}$.
 Note that  $\mathfrak{H}$ does not need to be  saturated  and  coincide with $\mathfrak{D}$. 

As equality of subschemes is a local property, we may check the equality locally. The proof is divided in two steps. 
\begin{enumerate}[{\it Step 1.}]
\item  For every Borel multi-index  ${\mathcal I}$ such that  $J(\mathcal I)\in \mathbb B$,  the ideal generated by 
$\mathfrak{h}$ defines ${\mathbf{H}}_{\mathcal{I}}$ as
closed subscheme 
of ${\mathbf{G}}_{\mathcal{I}}$.
\item For every (closed) point $I$ of $\GrassScheme{p}{N}$,  $\mathcal{Z}$ and  $\HilbScheme{p(t)}{n}$ coincide on a  neighborhood of $I$. 
\end{enumerate}

\noindent{\it Proof of Step 1.} We have to prove that for every  $k$-algebra $A$ and   ideal  $I$ in  $\underline{\mathbf{G}}_{\mathcal{I}}(A)$,  
$I$ is contained in $\underline{\mathbf{H}}_{\mathcal{I}}(A)$ if, and only if,  the polynomials in $\mathfrak{h}$ vanish when evaluated at $I$. 

Referring to Theorem \ref{thm:rifatto} and    Proposition \ref{solo2}, the vanishing at $I$ of the polynomials of $\mathfrak{h}_{1}$ 
is equivalent to  $\wedge^{q''+1}I^{(1)}=0$ and that of the  polynomials of $\mathfrak{h}_{2}$ 
 to $\wedge^{q''} I^{(1)}\wedge I''=0$.
 The equivalence  \emph{(\ref{it:rifatto_v})}$\Leftrightarrow$\emph{(\ref{it:rifatto_viii})}  of Theorem \ref{thm:rifatto} 
 and the definition of marked basis allow to conclude.

\medskip
  
\noindent{\it Proof of Step 2.}  Both ideals $\mathfrak{H}$ and $\mathfrak{D}$   are invariant under  the action of 
$\PGL$, $\mathfrak{H}$ by construction and $\mathfrak{D}$ because $\HilbScheme{p(t)}{n}$ is.

Due to the noetherianity of the ring of Pl\"ucker coordinates $k[\DD]$, we can choose     $ h_1, \dots, h_m\in \bigcup_{g \in \PGL} g\centerdot \mathfrak{h}$ 
that generate  $\mathfrak{H}$. If 
 $h_i\in g_i\centerdot\mathfrak{h}$, then  we get 
\[
 \left( g_1\centerdot\mathfrak{h}\cup \dots \cup  g_m\centerdot\mathfrak{h}\right)=\mathfrak{H}.
\] 
By the invariance of $ \mathfrak{H} $ under 
the action of $\PGL$,  we also get, for each $g\in \PGL$
\begin{equation*}
\left( gg_1\centerdot\mathfrak{h}\cup \dots
\cup g g_m\centerdot\mathfrak{h}\right)=g\centerdot\left( g_1\centerdot\mathfrak{h}\cup \dots \cup g_m\centerdot\mathfrak{h}\right)=g\centerdot\mathfrak{H} =
\mathfrak{H}. 
\end{equation*}

On the other hand,  if we restrict to the open subset  $\mathbf{G}_{\mathcal I,gg_1}\cap\dots\cap \mathbf{G}_{\mathcal I,gg_m}$, then by {\it Step 1}
and by the invariance of  $\mathfrak D$ under the action of $\PGL$,   we see that the ideal 
\[
\mathfrak D= \left( gg_1\centerdot\mathfrak{D}\cup \dots \cup g g_m\centerdot\mathfrak{D}\right)
\]
defines the same subscheme  as  $\mathfrak H=\left( gg_1\centerdot\mathfrak{h}\cup \dots \cup g g_m\centerdot\mathfrak{h}\right)$. 
Therefore, 
\[
\HilbScheme{p(t)}{n}\cap (\mathbf{G}_{\mathcal I,gg_1}\cap\dots\cap \mathbf{G}_{\mathcal I,gg_m})
=\mathcal Z\cap (\mathbf{G}_{\mathcal I,gg_1}\cap\dots\cap \mathbf{G}_{\mathcal I,gg_m}).
\]

It remains to  prove that for every $I\in \GrassScheme{p}{N}$, 
we can find suitable $g \in \PGL$ and   $J(\mathcal I)\in \mathbb B$, such that
  $I\in \mathbf{G}_{\mathcal I,gg_1}
\cap\dots\cap \mathbf{G}_{\mathcal I,gg_m}$.

By Proposition \ref{prop:GrBorelSubFunctors}, there are  $J({\mathcal I})\in \mathbb B$ and $\overline g$ such that $I\in {\mathbf{G}} _{{\mathcal{I}},\overline{g}}$.
 The orbit of $I$ under the action of $\PGL$ 
 is almost completely contained in
 ${\mathbf{G}}_{{\mathcal{I}},\overline{g}}$; let $U$ be an open subset of $\PGL$ such that
 $(g')^{-1}\centerdot I\in {\mathbf{G}}_{{\mathcal{I}},\overline{g}}$, i.e.~$I\in  {\mathbf{G}}_{{\mathcal{I}},g'\overline{g}}$.  
 Therefore, 
for a general $g\in \PGL$,  it holds   $g g_1\overline{g}^{-1}, \dots, g g_m\overline{g}^{-1}\in U$ and   $I\in \mathbf{G}_{\mathcal I,gg_1}\cap\dots\cap 
\mathbf{G}_{\mathcal I,gg_m}$ as wanted.
\end{proof}
 
 \bigskip

For sake of completeness we now show how our strategy also allows to mimic the construction of equations for $\HilbScheme{p(t)}{n}$ presented in the well-known papers by Iarrobino and Kleiman \cite{IarrobinoKleiman} and  by Haiman
and Sturmfels \cite{HaimSturm}.

\subsection{Equations of higher degree}\label{subsec:IK}

Let $A$ be a $k$-algebra and $I$ be  an ideal in $\underline{\mathbf{G}}_{{\mathcal{I}}}(A)$. Exploiting 
  Theorem \ref{universal},   we  obtain a set of generators for $I_{r+1}$ evaluating at $I$ the following set of polynomials
\begin{equation*}
x_0 \mathcal{B}^{(1)} \cup \cdots \cup x_n \mathcal{B}^{(1)}=\left\{x_i \delta^{(1)}_{\mathcal{K}}\ \big\vert\ i=0,\ldots,n,\ \mathcal{K}
 \in \mathcal{E}^{(1)} \right\} .
\end{equation*}
By Theorem \ref{thm:rifatto}\emph{(\ref{it:rifatto_v})}$\Leftrightarrow$\emph{(\ref{it:rifatto_vi})},  we know that 
 $I \in \underline{\mathbf{H}}_{{\mathcal{I}}}(A)$ if, and only if, 
$\wedge^{q(r+1)+1} I_{r+1}$ vanishes. 
The exterior power $\wedge^{q(r+1)+1} I_{r+1}$ is generated by all the possible exterior products of order $q(r+1)+1$ among the given set of generators 
of $I_{r+1}$. Therefore, 
the conditions  $I \in \underline{\mathbf{H}}_{\mathcal{I}}(A)$
is given by the vanishing at $I$ of the $x$-coefficients in the wedge
 products
\begin{equation*}
\bigwedge_{j=1}^{q(r+1)+1} x_{i_j} \delta^{(1)}_{\mathcal{K}_j},\qquad \forall \ 0 \leq i_1 \leq \dots \leq  i_{q(r+1)+1}\leq  n,\quad \forall\ 
\mathcal{K}_j \in \mathcal{E}^{(1)}.
\end{equation*}

The open subfunctors $\underline{\mathbf{G}}_{\mathcal{I}}$ cover the Grassmann functor and each $\underline{\mathbf{H}}_{\mathcal{I}}$ 
is representable, so that we can apply Proposition \ref{prop:openSubfunctors}. The natural transformations $\mathscr{H}_{\mathcal{I}}\colon
 \underline{\mathbf{H}}_{\mathcal{I}} \rightarrow \underline{\mathbf{G}}_{\mathcal{I}}$ are induced by closed embeddings of schemes, hence the same holds true
for $\mathscr{H}\colon\HilbFunctor{p(t)}{n} \rightarrow \GrassFunctor{p}{N}$.

\begin{theorem}[Iarrobino-Kleiman-like equations for the Hilbert scheme]\label{th:IarroKleiman}
The subscheme of $\GrassScheme{p}{N}$ representing the Hilbert functor $\HilbFunctor{p(t)}{n}$ can be defined by 
an ideal generated  by homogeneous elements  of degree $ q(r+1)+1$ in the ring $k[\DD]$ of the  Pl\"ucker coordinates.
\end{theorem}

The above equations of degree $q(r+1)+1$   coincides on each open subscheme  $\mathbf{G}_{\mathcal{I}}$ of the standard open 
cover of the Grassmannian with  the set of equations obtained  by Iarrobino and Kleiman in local coordinates.  We could also exploit this same 
argument  using the Borel open cover of $\GrassScheme{p}{N}$, instead of the standard one and obtain a different set of equations of the same degree.

\subsection{\texorpdfstring{Equations of degree $n+1$}{Equations of degree n+1}} \label{subsec:HS}  
As  pointed out by Haiman and Sturmfels, if  $I=(I_r)$ is generated by a set of polynomials $B$, then the matrix $M_{r+1}$ that represents the generators 
 $x_0 B \cup \cdots \cup x_n B$ of the module $I_{r+1}$  contains $n+1$ copies of the matrix $M_r$ corresponding to
 $B$. Hence,    some minors of $M_{r+1}$ are also
  minors of $M_r$   and every minor of $M_{r+1}$ can be obtained as the sum of products of at most $n+1$ minors of $M_r$.   
  
This observation suggests to expand $\wedge^{q(r+1)+1} I_{r+1}$ as done in  Remark \ref{rk:generat} and take the $\xx$-coefficients of
\begin{equation*}
  \bigwedge_{\begin{subarray}{c} 0 \leq i \leq n \\ m_i > 0\end{subarray}} x_i
\delta^{(m_i)}_{\mathcal{K}_i}  
  ,\qquad 	\forall \ \delta^{(m_i)}_{\mathcal{K}_i}  \in \mathcal{B}^{(m_i)}  \text{ s.t.~}    \sum m_i=m.  
\end{equation*}
\begin{theorem}[Bayer-Haiman-Sturmfels-like equations for the Hilbert scheme]\label{th:BayerHaimanSturmfels}
The subscheme of $\GrassScheme{p}{N}$ representing the Hilbert functor $\HilbFunctor{p(t)}{n}$ can be defined by an ideal generated
by homogeneous elements  of degree $\leq n+1$ in the ring $k[\DD]$ of the  Pl\"ucker coordinates.
\end{theorem}

In this case, if we use the standard open cover of $\GrassScheme{p}{N}$,  we obtain the same global equations given by Haiman and Sturmfels,   
while  using the Borel open cover we obtain a different set of equations with  maximum degree $n+1$.

\section{Examples: Hilbert schemes of points}\label{sec:example}

\subsection{\texorpdfstring{The Hilbert scheme $\HilbScheme{2}{2}$}{The Hilbert scheme of 2 points in the plane}}
The Gotzmann number of the Hilbert polynomial $p(t)=2$ is $r=2$, hence $N(r)=6$ and $p(r)=2$. We identify $H^0 \big(\mathcal O_{\PP^2}(2)\big)$ with $k^6$ by setting 
$\a_i=x^{\alpha(i)}$ where  $x^{\alpha(i)}$ is the $i$-th term in the sequence
 $(x_2^2,x_2x_1,x_1^2,x_2x_0,x_1x_0,x_0^2)$.
In this way we obtain the natural transformation of functors $\HilbFunctor{2}{2}
\rightarrow 
\GrassFunctor{2}{6}$. 

\medskip
\paragraph{\bf Standard open cover of $\GrassFunctor{2}{6}$}  There are $\binom{6}{2}=15$ open subfunctor $\underline{\mathbf{G}}_{\mathcal{I}}$ 
in the standard open cover, each corresponding to a $2$-multi-index
$\mathcal I\subset \{1,2,3,4,5,6\}$.  Not every   element of $\GrassFunctor{2}{6}(A)$ is  contained in one of them (not even the free ones),  if 
$A$ is not a field or even a local ring.

Let us consider for instance $A:=k[t]$ and 
\begin{equation}\label{eq:morph}
\pi\colon A^6 
\xrightarrow{\text{\tiny$\left(
\begin{array}{cccccc}
1-t & 0 & t^2 & 0 &0 & 0 \\
0 & 0 &1 &0 & 1+t & 1\end{array}\right)$}} 
 A^2.
\end{equation} 
This map is surjective, since  $ (1,0)=(1+t)\pi(\a_1)+\pi(\a_3)-\pi(\a_6)$ and $(0,1)=\pi(\a_6)$. 
Its  kernel  is the free $A$-module $L=\langle t^2\a_1+( t-1)\a_3-(t-1)\a_6,\ \a_2, \ \a_4, \ \a_5-(t+1)\a_6\rangle$.

\medskip

Thus, the quotient $Q:=A^6/\ker \pi$ is isomorphic to $A^2$ and $Q\in \GrassFunctor{2}{6}(A)$. Notice that the set of non-zero 
maximal minors $\{1-t,1-t^2,t^2+t^3,t^2\}$ 
of the matrix defining $\pi$ generates $A$, but none of them alone does, so that $Q$ does not belong to any 
$\underline{\mathbf{G}}_{\mathcal{I}}(A)$.

\medskip

On the other hand,  for the local $k$-algebra $A':=k[t]_{(t)} $, the $A'$-module   $Q' := Q \otimes_A A' $ is in
$ \underline{\mathbf{G}}_{\mathcal I}(A')$ for $\mathcal I=(1,3)$. 
The $\mathcal I$-marked set  of the $A'$-module $L'$ such that $Q' = A'^6/L'$ is
 \begin{equation*}
 \begin{array}{lll}
\b_2=  \a_2& & \b_5= \a_5 + \tfrac{t^3+t^2}{1-t}\a_1 - (1+t)\a_3 \\
\b_4= \a_4  && \b_6= \a_6 +\tfrac{t^2}{1-t}\a_1 - \a_3  \\
 \end{array}.
 \end{equation*}
The  Pl\"ucker coordinates of $L'$ (with $ \Delta_{13}(L') = 1$)  are given by Corollary \ref{prop:genDelta} 
{\normalsize
\[
\begin{array}{llllll}
 \Delta_{12}(L') = 0, && \Delta_{23}(L') = 0, &&& \Delta_{35}(L') = \frac{t^3+t^2}{1-t},\\ 
 \Delta_{13}(L') = 1, && \Delta_{24}(L') = 0, &&& \Delta_{36}(L') = \frac{t^2}{1-t},\\
 \Delta_{14}(L') = 0, && \Delta_{25}(L') = 0, &&& \Delta_{45}(L') = 0,\\
 \Delta_{15}(L') =( 1+t), && \Delta_{26}(L') = 0, &&& \Delta_{46}(L') = 0,\\
 \Delta_{16}(L') = 1, && \Delta_{34}(L') = 0, &&& \Delta_{56}(L') =0.\\
  \end{array}
\]
}

\medskip

\paragraph{\bf The generators $\mathcal{B}^{(m)}$} For $m=1$ there are 20 elements  in $\mathcal{B}^{(1)}$, since there are $\binom{6}{3}=20$
multi-indices $\mathcal K \in \mathcal{E}^{(1)}$. 
For instance for $\mathcal K \in \mathcal{E}_{13}^{(1)}$ we get 
\[
\begin{split}
\delta^{(1)}_{123}  &{}= \Delta_{23}\, \a_1 - \Delta_{13}\, \a_2 + \Delta_{12}\, \a_3  \\
\delta^{(1)}_{134}  &{}= \Delta_{34}\, \a_1  - \Delta_{14} \a_3 +\Delta_{13}\, \a_4  \\
\delta^{(1)}_{135}  &{}= \Delta_{35}\, \a_1 - \Delta_{15}\,  \a_3 + \Delta_{13}\, \a_5 \\
\delta^{(1)}_{136}  &{}= \Delta_{36}\, \a_1  - \Delta_{16} \, \a_3 + \Delta_{13}\,  \a_6\\
\end{split}
\]
and for   $\mathcal K =(3,5,6)$ we get 
\[
\delta^{(1)}_{356}  = \Delta_{56} \a_3-\Delta_{36} \a_5 + \Delta_{35} \a_6  .
\]
They are not independent. For instance there is the relation
\[  \Delta_{13}\, \delta^{(1)}_{356}  + \Delta_{36}\, \delta^{(1)}_{135}  -\Delta_{35} \delta^{(1)}_{136}
=(\Delta_{13}\,\Delta_{56}-\Delta_{15}\,\Delta_{36}+\Delta_{16}\,\Delta_{35}\,) \a_3=0 \]
(note that the expression in the round brackets is a Pl\"ucker relation).

For $m=2$, $\mathcal{B}^{(2)}$ contains  $\binom{6}{4} = 15$ elements. For instance,
\[
\delta^{(2)}_{1356}=\Delta_{13}\a_5 \wedge \a_6 - \Delta_{15} \a_3 \wedge \a_6 + \Delta_{16} \a_3 \wedge \a_5 + 
 \Delta_{35}\a_1 \wedge \a_6 - \Delta_{36}\a_1  \wedge \a_5 + \Delta_{56}\a_1\wedge \a_3.
\]
Finally,  $\mathcal{B}^{(3)}$  contains  $\binom{6}{5} = 6$ 
elements and $\mathcal{B}^{(4)}$ has a unique element.

\medskip

\paragraph{\bf Borel open cover of $\GrassFunctor{2}{6}$}

It is easy to check that there is only one Borel multi-index of two elements, namely $\mathcal I=(5,6)$. 

As the minor of the matrix \eqref{eq:morph} corresponding to the last two columns is identically zero, for every $\mathfrak{p} \in 
\Spec k[t]$, $Q_{\mathfrak{p}}$ 
is not contained in $\underline{\mathbf{G}}_{56}(k[t]_\mathfrak{p})$. We can apply Proposition \ref{prop:GrBorelSubFunctors} and determine 
for $\mathfrak{p} = (1-t)$ a change
of coordinates $g \in \PGL_{\QQ}(3)$ such that $Q_{\mathfrak{p}}$ is contained in $\underline{\mathbf{G}}_{56,g}(k(\mathfrak{p}))$.
Tensoring by the residue field $k(\mathfrak{p}) \simeq k$, we obtain the following surjective morphism of vector spaces
\[
k^6
\xrightarrow{\text{\tiny$\left(
\begin{array}{cccccc}
0 & 0 & 1 & 0 &0 & 0 \\
0 & 0 &1 &0 & 2 & 1\end{array}\right)$}} 
k^2
\]
whose kernel is the vector space $\langle x_2^2, x_2x_1, x_2x_0,x_1x_0 -x_0^2\rangle$. The generic initial
 ideal of the ideal $I = (x_2^2, x_2x_1, x_2x_0,x_1x_0 -x_0^2)$  is $J = (x_2^2,x_2x_1,x_1^2,x_2x_0)$. A change of coordinates
 $g$ such that $g \centerdot I = J$ is, for instance, the automorphism swapping $x_1$ and $x_0$. Indeed,
\[
g\centerdot (x_2^2, x_2x_1, x_2x_0,x_1x_0 -x_0^2) = (x_2^2, x_2x_0, x_2x_1,x_1^2 - \frac{1}{2}x_1x_0)
\] 
and 
\[
\widetilde{g} = \left(
\begin{array}{cccccc}
1 & 0 & 0 & 0 & 0 & 0 \\
0 & 0 & 0 & 1 & 0 & 0 \\
0 & 0 & 0 & 0 & 0 & 1 \\
0 & 1 & 0 & 0 & 0 & 0 \\
0 & 0 & 0 & 0 & 1 & 0 \\
0 & 0 & 1 & 0 & 0 & 0
\end{array}\right)
\quad
\text{so that}
\quad
\pi_{\mathfrak{p}} \circ \widetilde{g} \circ \Gamma_{56} = 
\left(
\begin{array}{cc}
0 & 1 \\ 2& 1
\end{array}\right)
\]
is surjective. Notice that this change of coordinates does not work for all localizations. Indeed
\[
\pi \circ \widetilde{g} \circ \Gamma_{56} = \left(
\begin{array}{cc}
0 & t^2 \\ t+1& 1
\end{array}\right)
\]
is not surjective in the localizations  $k[t]_{(t+1)}$ and $k[t]_{(t)}$ as the determinant is not invertible. 

\medskip

\paragraph{\bf New equations}

Let us finally show how to determine the equations of degree $2$ defining the scheme representing the Hilbert functor $\HilbFunctor{2}{2}$ 
in the Grassmannian
$\GrassScheme{2}{6}$. 

As $d=0$,   $I^{(1)}=x_0 I_2$ and its rank is equal to  $q(2) = q''=4$. Therefore, the 
set of equations $\mathfrak{h}_1$ of
\eqref{eq:part1} is empty.

The set of equations  $\mathfrak{h}_2$ ensures that  $\bigwedge^4 I^{(1)}\wedge I'' =0$ and contains 
the $\xx$-coefficients of  the products between $x_0 \delta^{(4)}_{123456}$ and each
element of $\mathcal{W} = \{ x_2\delta^{(1)}_{245} - x_1\delta^{(1)}_{145}, x_2\delta^{(1)}_{246} - x_1\delta^{(1)}_{146}, x_2\delta^{(1)}_{256} - x_1\delta^{(1)}_{156}, x_2\delta^{(1)}_{345} - x_1\delta^{(1)}_{245}, x_2\delta^{(1)}_{346} - x_1\delta^{(1)}_{246}, x_2\delta^{(1)}_{356} - x_1\delta^{(1)}_{256}, x_1 \delta^{(1)}_{456} - x_0\delta^{(1)}_{256}, x_2 \delta^{(1)}_{456} - x_0\delta^{(1)}_{156}, x_1 \delta^{(1)}_{456} + x_0\delta^{(1)}_{346}, x_2 \delta^{(1)}_{456} + x_0\delta^{(1)}_{246}\}$. We notice that this set gives redundant relations, indeed for instance 
\small
\[
\begin{split}
x_0 \delta^{(4)}_{123456} &{}\wedge \left(x_1 \delta^{(1)}_{456} - x_0\delta^{(1)}_{256}\right) = x_0 \delta^{(4)}_{123456} \wedge x_1 \delta^{(1)}_{456} - x_0 \left(\delta^{(4)}_{123456} \wedge \delta^{(1)}_{256}\right) = x_0 \delta^{(4)}_{123456} \wedge x_1 \delta^{(1)}_{456} = {}\\
&{}= x_0 \delta^{(4)}_{123456} \wedge x_1 \delta^{(1)}_{456} + x_0 \left(\delta^{(4)}_{123456} \wedge \delta^{(1)}_{346}\right) = x_0 \delta^{(4)}_{123456} \wedge \left(x_1 \delta^{(1)}_{456} + x_0\delta^{(1)}_{346}\right) 
\end{split}
\]
\normalsize
as each exterior product of order greater than 4 vanishes (in fact the $x$-coefficients we obtain from such products of order 5 are contained in the ideal generated by the Pl\"ucker relations). Hence, in order to determine the equations of $\mathfrak{h}_2$ we consider the set of polynomials $\widehat{\mathcal{W}} = \{ x_2\delta^{(1)}_{245} - x_1\delta^{(1)}_{145}, x_2\delta^{(1)}_{246} - x_1\delta^{(1)}_{146}, x_2\delta^{(1)}_{256} - x_1\delta^{(1)}_{156}, x_2\delta^{(1)}_{345} - x_1\delta^{(1)}_{245}, x_2\delta^{(1)}_{346} - x_1\delta^{(1)}_{246}, x_2\delta^{(1)}_{356} - x_1\delta^{(1)}_{256}, x_1 \delta^{(1)}_{456}, x_2 \delta^{(1)}_{456}\}$.
We get 48 equations which are reduced to 30 by Pl\"ucker relations. 

To obtain the equation defining $\HilbScheme{2}{2} \subset \GrassScheme{2}{6}$, we need to determine the orbit of these polynomials with respect to the action of $\PGL_{\QQ}(3)$. 
However, in this special case, we discover that the ideal generated by the Pl\"ucker relations and by $\mathfrak{h}_2$ is already $\PGL_{\QQ}(3)$ invariant, i.e.~the equations in $\mathfrak{h}_2$ define the Hilbert scheme. The Pl\"ucker relations and the following $30$ equations define $\HilbScheme{2}{2}$ as subscheme of $\PP^{14}$:
\small
\[
\begin{split}
& \Delta_{13} \Delta_{14}-\Delta_{12} \Delta_{24}-\Delta_{12} \Delta_{15},\
\Delta_{13} \Delta_{24}-\Delta_{12} \Delta_{34}-\Delta_{12} \Delta_{25},\
\Delta_{23} \Delta_{24}-\Delta_{12} \Delta_{35},\\
& \Delta_{14} \Delta_{24}-\Delta_{14} \Delta_{15}+\Delta_{12} \Delta_{16},\
\Delta_{24}^2-\Delta_{14} \Delta_{25}+\Delta_{12} \Delta_{26},\
\Delta_{23} \Delta_{34}+\Delta_{23} \Delta_{25}-\Delta_{13} \Delta_{35},\\
& \Delta_{14} \Delta_{34}-\Delta_{14} \Delta_{25}+\Delta_{12} \Delta_{45}+\Delta_{12} \Delta_{26},\
\Delta_{24} \Delta_{34}-\Delta_{14} \Delta_{35}+\Delta_{12} \Delta_{36},\\
& \Delta_{34}^2-\Delta_{15} \Delta_{35}-\Delta_{23} \Delta_{45}+\Delta_{13} \Delta_{36},\
\Delta_{15}^2-\Delta_{14} \Delta_{25}-\Delta_{13} \Delta_{16}+\Delta_{12} \Delta_{26},\
\Delta_{24} \Delta_{25}-\Delta_{14} \Delta_{35},\\
&\Delta_{15} \Delta_{25}-\Delta_{14} \Delta_{35}-\Delta_{13} \Delta_{26}+\Delta_{12} \Delta_{36},\
\Delta_{25}^2-\Delta_{15} \Delta_{35}-\Delta_{23} \Delta_{45},\\
& \Delta_{24} \Delta_{35}-\Delta_{15} \Delta_{35}-\Delta_{23} \Delta_{45}+\Delta_{23} \Delta_{26},\
\Delta_{34} \Delta_{35}-\Delta_{25} \Delta_{35}+\Delta_{23} \Delta_{36},\
\Delta_{14} \Delta_{45}-\Delta_{12} \Delta_{46},\\
& \Delta_{24} \Delta_{45}-\Delta_{12} \Delta_{56},\
\Delta_{34} \Delta_{45}+\Delta_{23} \Delta_{46}-\Delta_{13} \Delta_{56},\
\Delta_{15} \Delta_{45}-\Delta_{13} \Delta_{46}+\Delta_{12} \Delta_{56},\\
& \Delta_{25} \Delta_{45}-\Delta_{23} \Delta_{46},\
\Delta_{35} \Delta_{45}-\Delta_{23} \Delta_{56},\
\Delta_{45}^2-\Delta_{25} \Delta_{46}+\Delta_{15} \Delta_{56},\\
& \Delta_{24} \Delta_{26}-\Delta_{14} \Delta_{36}+\Delta_{12} \Delta_{56},\
\Delta_{25} \Delta_{26}-\Delta_{15} \Delta_{36}-\Delta_{23} \Delta_{46}+\Delta_{13} \Delta_{56},\\
& \Delta_{26}^2-\Delta_{16} \Delta_{36}-\Delta_{25} \Delta_{46}+\Delta_{15} \Delta_{56},\
\Delta_{24} \Delta_{46}-\Delta_{14} \Delta_{56},\
\Delta_{34} \Delta_{46}+\Delta_{25} \Delta_{46}-\Delta_{24} \Delta_{56}-\Delta_{15} \Delta_{56},\\
&\Delta_{35} \Delta_{46}-\Delta_{25} \Delta_{56},\
\Delta_{45} \Delta_{46}-\Delta_{26} \Delta_{46}+\Delta_{16} \Delta_{56},\
\Delta_{36} \Delta_{46}-\Delta_{45} \Delta_{56}-\Delta_{26} \Delta_{56}.\\
\end{split}
\]
\normalsize
Furthermore, we check that the ideal they generate is saturated, then it is \emph{the saturated ideal} of $\HilbScheme{2}{2}$. Its Hilbert polynomial is $\frac{21}{4!}t^4+\frac{15}{4}t^3+\frac{45}{8}t^2+\frac{15}{4}t+1$, hence $\HilbScheme{2}{2} \subset \PP^{14}$ is a subscheme of dimension $4$ (as expected) and  degree $21$,  as already proved  in \cite{HaimSturm,BrodskySturmfels}. A different set of quadratic equations defining this Hilbert scheme can be obtained also using border bases and commutation relations of multiplicative matrices (see \cite{ABM}).

\smallskip

\paragraph{\bf Iarrobino-Kleiman equations}

Let us now see how to compute the Iarrobino-Kleiman equations for $\HilbScheme{2}{2}$.  The universal element of   
$\GrassScheme{2}{6}$  is generated by  $\mathcal{B}^{(1)}$. In order to compute   $\wedge^{q(r+1)+1} I_{r+1}=\wedge^{9} I_{3}$ 
we use the set of generators  $x_0\mathcal{B}^{(1)} \cup x_1\mathcal{B}^{(1)}
\cup x_2\mathcal{B}^{(1)}$ of $I_3$. The $\xx$-coefficients of any exterior product of order $9$  are 
expressions  of degree $9$ in the Pl\"ucker coordinates. Their union  defines
$\HilbScheme{2}{2} \subset \GrassScheme{2}{6}$. 

For instance, considering the $9$ elements $x_2 \delta^{(1)}_{126}$, 
$x_2 \delta^{(1)}_{156}$, $x_2 \delta^{(1)}_{234}$, $x_2 \delta^{(1)}_{356}$, $x_1 \delta^{(1)}_{123}$, $x_1 \delta^{(1)}_{345}$, $x_0 
\delta^{(1)}_{146}$, $x_0 
\delta^{(1)}_{234}$,  $x_0 \delta^{(1)}_{456}$, the $x$-coefficients of their exterior product are the maximal minors of the following matrix
 \footnotesize
 \begin{equation*}
 \begin{split}
 & \quad \begin{array}{>{$}p{0.7cm}<{$} >{$}p{0.8cm}<{$} >{$}p{0.8cm}<{$} >{$}p{0.7cm}<{$} >{$}p{0.8cm}<{$} >{$}p{1cm}<{$}>{$}p{0.8cm}<{$}>{$}p{0.8cm}<{$} >{$}p{0.8cm}<{$} >{$}p{0.7cm}<{$}}
 x_2^3 & x_2^2 x_1 & x_2 x_1^2 & x_1^3 & x_2^2 x_0 & x_2 x_1 x_0 & x_1^2 x_0 & x_2 x_0^2 & x_1 x_0^2 & x_0^3\\
 \end{array}\\
 \begin{array}{r}	
 x_2 \delta^{(1)}_{126} \\
 x_2 \delta^{(1)}_{156} \\
 x_2 \delta^{(1)}_{234} \\
 x_2 \delta^{(1)}_{356} \\
 x_1 \delta^{(1)}_{123} \\
 x_1 \delta^{(1)}_{345} \\
 x_0 \delta^{(1)}_{146} \\
 x_0 \delta^{(1)}_{234} \\
 x_0 \delta^{(1)}_{456}
 \end{array}
 & \left(
  \begin{array}{>{$}p{0.7cm}<{$} >{$}p{0.8cm}<{$} >{$}p{0.8cm}<{$} >{$}p{0.7cm}<{$} >{$}p{0.8cm}<{$} >{$}p{1cm}<{$} >{$}p{0.8cm}<{$}>{$}p{0.8cm}<{$} >{$}p{0.8cm}<{$} >{$}p{0.7cm}<{$}}
 \Delta_{26}& {-\Delta_{16}} & 0& 0& 0& 0& 0& \Delta_{12}& 0\phantom{\delta^{(}_1} & 0\\
  \Delta_{56} & 0& 0& 0& 0& {-\Delta_{16}}& 0& \Delta_{15}& 0\phantom{\delta^{(}_1} & 0 \\
 0\phantom{\delta^{(}_1}& \Delta_{34}& {-\Delta_{24}}& 0& \Delta_{23}& 0& 0& 0& 0& 0\\
  0\phantom{\delta^{(}_1}& 0& \Delta_{56}& 0& 0& {-\Delta_{36}}& 0& \Delta_{35}& 0& 0\\
  0\phantom{\delta^{(}_1}& \Delta_{23}& {-\Delta_{13}}& \Delta_{12}& 0& 0& 0& 0& 0& 0\\
 0\phantom{\delta^{(}_1}& 0& 0& \Delta_{45}& 0& {-\Delta_{35}}& \Delta_{34}& 0& 0& 0\\
  0\phantom{\delta^{(}_1}& 0& 0& 0& \Delta_{46}& 0& 0& {-\Delta_{16}}& 0& \Delta_{14}\\
  0\phantom{\delta^{(}_1}& 0& 0& 0& 0& \Delta_{34}& {-\Delta_{24}}& \Delta_{23}& 0& 0\\
  0\phantom{\delta^{(}_1}& 0& 0& 0& 0& 0& 0& \Delta_{56}& {-\Delta_{46}}& \Delta_{45}\end{array}\right).\end{split}
 \end{equation*}
 \normalsize

\medskip

\paragraph{\bf Bayer-Haiman-Sturmfels equations}

In order to lower the degree of the equations, we can impose the vanishing of the exterior power $\wedge^9 I_3$ by considering $\wedge^9 I_3$ 
generated by all
exterior products $x_0 \delta^{(m_0)}_{\mathcal{J}_0} \wedge x_1 \delta^{(m_1)}_{\mathcal{J}_1} \wedge x_2 
\delta^{(m_2)}_{\mathcal{J}_2}$ for $(m_0,m_1,m_2)$ such that $m_0 + m_1 + m_2 = 9$ and $0 \leq m_0,m_1,m_2 \leq 4$. For $m_0=
3$, $m_1=2$, $m_2=4$,  we get for instance  $x_0 \delta^{(3)}_{23456} \wedge x_1 \delta^{(2)}_{1346} \wedge x_2 \delta^{(4)}_{123456} $, where

\scriptsize
\[
 \begin{split}
 x_0 \delta^{(3)}_{23456} ={}& \Delta_{56}\, x_2x_1x_0 \wedge x_1^2x_0 \wedge x_2x_0^2 - \Delta_{46}\, x_2x_1x_0 \wedge x_1^2x_0 \wedge x_1x_0^2 
+ \Delta_{45}\, x_2x_1x_0\wedge x_1^2x_0 \wedge x_0^3 + {}\\ &\Delta_{36}\, x_2x_1x_0 \wedge x_2x_0^2 \wedge x_1x_0^2 - \Delta_{35}\, x_2x_1x_0 
\wedge x_2x_0^2 \wedge x_0^3 + \Delta_{34}\, x_2x_1x_0 \wedge x_1x_0^2 \wedge x_0^3 -{}\\ & \Delta_{26}\, x_1^2x_0 \wedge x_2x_0^2 \wedge x_1x_0^2
 + \Delta_{25}\, x_1^2x_0 \wedge x_1x_0^2 \wedge x_0^3 - \Delta_{24}\, x_1^2x_0 \wedge x_1x_0^2 \wedge x_0^3 +\Delta_{23}\, x_2x_0^2 
\wedge x_1x_0^2 \wedge x_0^3,\\
  x_1 \delta^{(2)}_{1346} = {}& \Delta_{46}\, x_2^2 x_1 \wedge x_1^3 - \Delta_{36}\, x_2^2 x_1 \wedge x_2x_1x_0 + \Delta_{34}\, x_2^2 x_1
  \wedge x_1x_0^2 +
  \Delta_{16}\, x_1^3 \wedge x_2x_1x_0 - \Delta_{14}\, x_1^3 \wedge x_1x_0^2 + {}\\ &\Delta_{13}\, x_2x_1x_0 \wedge x_1x_0^2,\\
  x_2 \delta^{(4)}_{123456}  ={}& \Delta_{56}\, x_2^3 \wedge x_2^2 x_1 \wedge x_2x_1^2 \wedge x_2^2 x_0 - \Delta_{46}\, x_2^3 \wedge x_2^2 x_1 
\wedge x_2x_1^2 \wedge x_2x_1x_0 + \Delta_{45}\, x_2^3 \wedge x_2^2 x_1 \wedge x_2x_1^2 \wedge x_2x_0^2 + {}\\
& \Delta_{36}\, x_2^3 \wedge x_2^2 x_1 \wedge x_2^2 x_0 \wedge 
x_2x_1x_0 - \Delta_{35}\, x_2^3 \wedge x_2^2 x_1 \wedge x_2^2 x_0 \wedge x_2x_0^2 + \Delta_{34}\, x_2^3 \wedge x_2^2 x_1 \wedge
 x_2x_1x_0 \wedge x_2x_0^2 - {}\\
  & \Delta_{26}\, x_2^3 \wedge x_2x_1^2 \wedge x_2^2 x_0 \wedge x_2x_1x_0 + \Delta_{25}\, x_2^3 \wedge x_2x_1^2 \wedge x_2^2 x_0 \wedge
 x_2x_0^2 -  \Delta_{24}\, x_2^3 \wedge x_2x_1^2 \wedge x_2x_1x_0 \wedge x_2x_0^2 + {}\\ &\Delta_{23}\, x_2^3 \wedge x_2^2 x_0 \wedge
 x_2x_1x_0 \wedge x_2x_0^2 + \Delta_{16}\, x_2^2 x_1 \wedge x_2x_1^2 \wedge x_2^2 x_0 \wedge x_2x_1x_0 -\Delta_{15}\, x_2^2 x_1 
\wedge x_2x_1^2 \wedge x_2^2 x_0 \wedge x_2x_0^2 + {}\\
  & \Delta_{14}\, x_2^2 x_1 \wedge x_2x_1^2 \wedge x_2x_1x_0 \wedge x_2x_0^2 - \Delta_{13}\, x_2^2 x_1 \wedge x_2^2 x_0 \wedge
 x_2x_1x_0
 \wedge x_2x_0^2 + \Delta_{12}\, x_2x_1^2 \wedge x_2^2 x_0 \wedge x_1^2x_0 \wedge x_2x_0^2.\\
 \end{split}
\]
\normalsize
Its $\xx$-coefficients are the following polynomials  of degree 3 in the Pl\"ucker coordinates:
\small
\[
 \begin{split}
& \Delta_{26}^2 \Delta_{46}- \Delta_{25} \Delta_{46}^2- \Delta_{16}  \Delta_{26} \Delta_{56}+ \Delta_{14} \Delta_{56}^2,\  \Delta_{25}  \Delta_{26} \Delta_{46}- \Delta_{25} \Delta_{45} \Delta_{46}- \Delta_{16}  \Delta_{25} \Delta_{56},\\
& \Delta_{24}  \Delta_{26} \Delta_{46}- \Delta_{16}  \Delta_{24} \Delta_{56}- \Delta_{14} \Delta_{45} \Delta_{56},\ \Delta_{23}  \Delta_{26} \Delta_{46}+ \Delta_{25}  \Delta_{34} \Delta_{46}- \Delta_{16}  \Delta_{23} \Delta_{56}- \Delta_{14} \Delta_{35} \Delta_{56},\\
& \Delta_{25}  \Delta_{26}  \Delta_{34}- \Delta_{24}  \Delta_{25} \Delta_{36}- \Delta_{25}  \Delta_{34} \Delta_{45}+ \Delta_{13} \Delta_{25} \Delta_{56},\ \Delta_{16}  \Delta_{24} \Delta_{45}+ \Delta_{14} \Delta_{45}^2+ \Delta_{24}^2 \Delta_{46}- \Delta_{14}  \Delta_{25} \Delta_{46},\\
& \Delta_{24}  \Delta_{25} \Delta_{46}- \Delta_{14}  \Delta_{25} \Delta_{56},\ \Delta_{16}  \Delta_{24} \Delta_{35}- \Delta_{14}  \Delta_{25} \Delta_{36}+ \Delta_{14} \Delta_{35} \Delta_{45}+ \Delta_{23} \Delta_{24} \Delta_{46},\\
& \Delta_{16}  \Delta_{24}  \Delta_{25}- \Delta_{14}  \Delta_{25}  \Delta_{26}+ \Delta_{14}  \Delta_{25} \Delta_{45},\ \Delta_{15}  \Delta_{16}  \Delta_{24}- \Delta_{14}  \Delta_{16}  \Delta_{25}+ \Delta_{14}  \Delta_{15} \Delta_{45}- \Delta_{12}
  \Delta_{24} \Delta_{46}.
 \end{split}
\]
\normalsize

\medskip

In Table \ref{tab:compare}, there is a comparison between the number of  generators of the ideal defining the Hilbert scheme obtained according to the three different strategies.


\subsection{\texorpdfstring{The Hilbert scheme  $\HilbScheme{2}{3}$ }{The Hilbert scheme of 2 points in the space}}
The Hilbert scheme of 2 points in the projective space $\PP^3$ can be constructed as subscheme of the Grassmannian $\GrassScheme{2}{10} \subset \PP^{44}$. The set $\mathfrak{h}_1$ is empty (this happens for every Hilbert scheme of points) and the set $\mathfrak{h}_2$ contains $600$ equations of degree 2 that can be reduced to 330 modulo the Pl\"ucker relations. The ideal generated by $\mathfrak{h}_2$ and by the Pl\"ucker relations is not $\PGL_{\QQ}(4)$-invariant.
To obtain the equations defining $\HilbScheme{2}{3} \subset \GrassScheme{2}{10}$, we need to determine the orbits of these polynomials with respect to the action of $\PGL_{\QQ}(4)$. 
    From a computational point of view, we   consider a random element of $\PGL_{\QQ}(4)$, apply to our set of equations the induced automorphism  on the ring of Pl\"ucker coordinates of  $\GrassScheme{2}{10}$, add the  new equations to the previous set and repeat this process until the generated ideal stabilizes. 
    
The ideal we obtain is again saturated and its Hilbert polynomial is 
\[
\tfrac{370}{6!}t^6+\tfrac{83}{24}t^5+\tfrac{86}{9}t^4+\tfrac{335}{24}t^3+\tfrac{823}{72}t^2+\tfrac{61}{12}t+1
\]  
so that $\HilbScheme{2}{3}$ turns out to be a subscheme of $\PP^{44}$ of dimension $6$ and degree $370$,  defined by   $570$ quadratic equations ($210$ of them are Pl\"ucker relations). 

\subsection{\texorpdfstring{The Hilbert scheme  $\HilbScheme{2}{4}$ }{The Hilbert scheme of 2 points in the 4-space}}
The Hilbert scheme of 2 points in $\PP^4$ is constructed as subscheme of the Grassmannian $\GrassScheme{3}{15} \subset \PP^{104}$. From the computational point of view, the hardest part is the computation of the orbit of the equations $g\centerdot \mathfrak{h}$ for a given change of coordinates $g \in \PGL_{\QQ}(5)$. A first trick is to start considering simple changes of coordinates, for instance change of sign of a variable ($x_i \to -x_i$), swap of two variables ($x_i \leftrightarrow x_j$) and sum of two variables ($x_i \to x_i + x_j$). These changes of coordinates are easier to compute and bring us closer to the $\PGL_{\QQ}(5)$-invariant ideal of the Hilbert scheme, but in general they are not sufficient. In this case, a generic (random) change of coordinates $g \in \PGL_{\QQ}(5)$ induces a change of coordinates of $\GrassScheme{3}{15}$ described by a dense $105\times105$ matrix, so that computing the action of $g$ on a monomial of degree 2 in the Pl\"ucker coordinates requires more than 10000 multiplications, as each variable is replaced by a linear form with 105 terms. Therefore, it would be better to avoid redundancy in the equations of $\mathfrak{h}$. It is possible to reduce the redundancy replacing the set $\mathcal{B}^{(m_i)}$ with the union $\bigcup \mathcal{B}^{(m_i)}_{\mathcal{I}}$, with $\mathcal{I}$ a Borel multi-index, in the definition of equations \eqref{eq:part1} and \eqref{eq:part2}. In the case of 2 points, there is a unique Borel multi-index and for instance in the case of $\PP^4$, the set $\mathcal{B}^{(1)}$ contains $\binom{15}{3}$ polynomials while $\mathcal{B}^{(1)}_{14,15}$ has only $13$ polynomials. Applying these two tricks, we are able to compute the equations of $\HilbScheme{2}{4} \subset \PP^{104}$. The equations contained in $\mathfrak{h}$ obtained considering only the Borel multi-index are 480 and besides 24 simple changes of coordinates we need 3 random changes of coordinates to obtain the $\PGL_{\QQ}(5)$-invariant ideal. Finally, the Hilbert scheme turns out to be a subscheme defined by 3575 quadratic equations (1365 of them are Pl\"ucker relations) with Hilbert polynomial
\[
\tfrac{6125}{8!}t^8 + \tfrac{452}{288}t^7 + \tfrac{4027}{576}t^6 + \tfrac{635}{36}t^5 + \tfrac{31703}{1152}t^4 + \tfrac{7849}{288}t^3 + \tfrac{4849}{288} t^2 + \tfrac{145}{24}t + 1,
\]
i.e.~$\HilbScheme{2}{4}$ is a subscheme of $\PP^{104}$ of dimension 8 and degree 6125.

\subsection{\texorpdfstring{The Hilbert scheme  $\HilbScheme{3}{2}$ }{The Hilbert scheme of 3 points in the plane}}

The Hilbert scheme $\HilbScheme{3}{2}$ can be defined as subscheme of the Grassmannian $\GrassScheme{3}{10} \subset \PP^{119}$. There are two Borel-fixed ideals defining 3 points in the plane: $(x_2,x_1^3)$ and $(x_2^2,x_2x_1,x_1^3)$, so that in this case we can restrict to Borel multi-indices considering the elements of $\mathcal{B}^{(1)}_{7,9,10}$ and $\mathcal{B}^{(1)}_{8,9,10}$. In this way, the set $\mathfrak{h}$ contains $720$ equations and we obtain a $\PGL_{\QQ}(3)$-invariant ideal applying 10 changes of coordinates (8 special and 2 random). The ideal defining the Hilbert scheme is generated by 5425 quadratic equations (2310 Pl\"ucker relations) and $\HilbScheme{3}{2}$ is a subscheme of $\PP^{119}$ of dimension 6 and degree 3309, as its Hilbert polynomial is
\[
\tfrac{3309}{6!}t^6 + \tfrac{1557}{80}t^5 + \tfrac{553}{16}t^4 + \tfrac{543}{16}t^3 + \tfrac{2381}{120}t^2 + \tfrac{33}{5}t+1.
\]

\begin{table}[!ht]
\begin{center}
\begin{tikzpicture}[scale=0.84]
\node at (-1.9,0.1) [] {\footnotesize $\HilbScheme{2}{2} \subset \GrassScheme{2}{6} \subset \PP^{14}$};
\node at (-1.9,-0.4) [] {\tiny (15 Pl\"ucker relations)};
\node at (-1.9,-1.9) [] {\footnotesize $\HilbScheme{2}{3}\subset \GrassScheme{2}{10} \subset \PP^{44}$};
\node at (-1.9,-2.4) [] {\tiny (210 Pl\"ucker relations)};
\node at (-1.9,-3.9) [] {\footnotesize $\HilbScheme{2}{4}\subset \GrassScheme{2}{15} \subset \PP^{104}$};
\node at (-1.9,-4.4) [] {\tiny (1365 Pl\"ucker relations)};
\node at (-1.9,-5.9) [] {\footnotesize $\HilbScheme{3}{2}\subset \GrassScheme{3}{10} \subset \PP^{119}$};
\node at (-1.9,-6.4) [] {\tiny (2310 Pl\"ucker relations)};

\draw [thick] (0,2) -- (0,-7);
\draw [thick] (-3.8,1) -- (-3.8,-7);
\draw [thick] (5,2) -- (5,-7);
\draw [thick] (10,2) -- (10,-7);
\draw [thick] (15,2) -- (15,-7);
\draw [thick] (-3.8,1) -- (15,1);
\draw [thick] (0,2) -- (15,2);
\draw [thick] (-3.8,-1) -- (15,-1);
\draw [thick] (-3.8,-3) -- (15,-3);
\draw [thick] (-3.8,-5) -- (15,-5);
\draw [thick] (-3.8,-7) -- (15,-7);

\node at (2.5,1.5) [] {New equations};
\node at (12.5,1.5) [] {I-K equations};
\node at (7.5,1.5) [] {B-H-S equations};

\node at (2.5,0) [] {\tiny $\begin{array}{l} \text{Degree of equations: } 2\\ \dim \left(I_{\HilbScheme{2}{2}}\right)_2 = 45\\ \text{Number of equations: } \vert\mathfrak{h}\vert = 24\\ \text{Changes of coordinates: } 8+0\end{array}$};

\node at (2.5,-2) [] {\tiny $\begin{array}{l} \text{Degree of equations: } 2\\ \dim \left(I_{\HilbScheme{2}{3}}\right)_2 = 570\\ \text{Number of equations: } \vert\mathfrak{h}\vert = 140\\ \text{Changes of coordinates: } 15+1\end{array}$};

\node at (2.5,-4) [] {\tiny $\begin{array}{l} \text{Degree of equations: } 2\\ \dim \left(I_{\HilbScheme{2}{4}}\right)_2 = 3575\\ \text{Number of equations: } \vert\mathfrak{h}\vert = 480\\ \text{Changes of coordinates: } 24+3\end{array}$};

\node at (2.5,-6) [] {\tiny $\begin{array}{l} \text{Degree of equations: } 2\\ \dim \left(I_{\HilbScheme{3}{2}}\right)_2 = 5425\\ \text{Number of equations: } \vert\mathfrak{h}\vert = 720\\ \text{Changes of coordinates: } 8+2\end{array}$};

\node at (7.5,0) [] {\tiny $\begin{array}{l} \text{Degree of equations: } 3\\ \dim \left(I_{\HilbScheme{2}{2}}\right)_3 = 445 \\ \text{Number of equations: } \sim 8160\end{array}$};

\node at (7.5,-2) [] {\tiny $\begin{array}{l} \text{Degree of equations:} \leqslant 4\\ \dim \left(I_{\HilbScheme{2}{3}}\right)_4 = 185390\\ \text{Number of equations:} \sim 2{\cdot}10^{11} \end{array}$};

\node at (7.5,-4) [] {\tiny $\begin{array}{l} \text{Degree of equations:} \leqslant 5\\ \dim \left(I_{\HilbScheme{2}{4}}\right)_5 = 116461170\\ \text{Number of equations:} \sim 4{\cdot}10^{22} \end{array}$};

\node at (7.5,-6) [] {\tiny $\begin{array}{l} \text{Degree of equations:} \leqslant 3\\ \dim \left(I_{\HilbScheme{3}{2}}\right)_3 = 283245\\ \text{Number of equations:} \sim 6{\cdot}10^8 \end{array}$};

\node at (12.5,0) [] {\tiny $\begin{array}{l} \text{Degree of equations: } 9\\ \dim \left(I_{\HilbScheme{2}{2}}\right)_9 = 808225\\ \text{Number of equations:}\sim 9{\cdot}10^{10}\end{array}$};

\node at (12.5,-2) [] {\tiny $\begin{array}{l} \text{Degree of equations: } 19\\ \dim \left(I_{\HilbScheme{2}{3}}\right)_{19} \sim 6{\cdot}10^{15} \\ \text{Number of equations:} \sim 9{\cdot}10^{34} \end{array}$};

\node at (12.5,-4) [] {\tiny $\begin{array}{l} \text{Degree of equations: } 34\\ \dim \left(I_{\HilbScheme{2}{4}}\right)_{34} \sim 2{\cdot}10^{32} \\ \text{Number of equations:} \sim 10^{77} \end{array}$};

\node at (12.5,-6) [] {\tiny $\begin{array}{l} \text{Degree of equations: } 13\\ \dim \left(I_{\HilbScheme{3}{2}}\right)_{13} \sim 3{\cdot}10^{17} \\ \text{Number of equations:} \sim 3{\cdot}10^{28} \end{array}$};
\end{tikzpicture}
\end{center}
\caption{\label{tab:compare} A comparison among the characteristics of the different sets of equations defining the Hilbert schemes discussed in Section \ref{sec:example}. The set $\mathfrak{h}$ of new equations in this table contains the equations obtained considering only Borel multi-indices. In order to determine the $\PGL_{\QQ}(n+1)$-invariant ideal, we have always applied $(n+1)^2-1$ (the first summand) special changes of coordinates and the second summand corresponds to the needed random changes of coordinates (see \href{http://tinyurl.com/EquationsHilbPoints-m2}{\texttt{tinyurl.com/EquationsHilbPoints-m2}} for the explicit computations). The values of the Hilbert function of the ideals defining the Hilbert schemes have been computed from the ideal generated by the new equations. Notice the large redundancy of Bayer-Haiman-Sturmfels equations and Iarrobino-Kleiman equations which do not take into account the symmetries of the Hilbert scheme.}
\end{table}


\end{document}